%% file: word_basis-inline.tex
\def\filename{{word\_basis.tex}} 
\title{Pascal's triangle and word bases \\   for  blob algebra ideals} 
\author{P P Martin \\ \myaddress}\date{}
\documentclass[12pt]{article} 
\usepackage{epic} 
\usepackage{eepic} 
\usepackage{graphicx}
\usepackage{latexsym}
\usepackage{amssymb}
\input martinew
\newcommand{\beqa}{\begin{eqnarray}}%
\newcommand{\eeqa}{\end{eqnarray}}%

\newcommand{\U}{U}
\newcommand{\out}{\! \uparrow\!}
\newcommand{\back}{\downarrow}
\newcommand{\E}[1]{E_{#1}}
\newcommand{\Er}[1]{E'_{#1}}
\newcommand{\desc}{\setminus}
\newcommand{\ddesc}{\setminus\!\!\setminus}
\newcommand{\Powerset}{{\mathcal P}}%

\newcommand{\vardelta}{\delta_e}%

\newcommand{\ym}{{m}}

\begin{document} \maketitle 
 \newcommand{\ignoreifnotdraft}[1]{\ignore{#1}}
\ignoreifnotdraft{
\pagestyle{myheadings} \markboth{Draft}{\today}
}

\ignoreifnotdraft{
\noindent
\begin{tabular}{l} 
  {\tiny \filename (Draft)} \hspace{4.4in} Jan 1997 \\ \hline 
\end{tabular}
}
%
%
\ignoreifnotdraft{ \newpage } 

\section{Introduction}

The Temperley--Lieb algebra \cite{TemperleyLieb71} 
plays a role in many different branches of Mathematics and Physics
(see for example \cite{Baxter82,GoodmandelaHarpeJones89,Kauffman91,%
Martin91,Henkel99}). 
The blob algebra \cite{MartinSaleur94a} is a generalisation which preserves
this versatility   
(see for example 
\cite{MartinSaleur93,MartinSaleur94a,MartinWoodcock2000,%
MartinWoodcock03,MartinRyom02,deGier02,%
NicholsRittenbergdeGier05,
DoikouMartin03,DeGierPyatov03,GrahamLehrer03}). 

In particular, 
the blob algebra is used to study (affine) Hecke
algebra representation theory \cite{MartinWoodcock03,GrahamLehrer03},
and boundary integrable lattice models 
\cite{DoikouMartin03,DeGierPyatov03}.  
In these settings it is useful to have a
rather explicit understanding of the
relationship between the well-known `diagram' bases for 
`standard' blob modules 
\cite{MartinSaleur94a}  
and bases of words in abstract generators. 
This connection is sketched, in principle, both in
\cite{MartinSaleur94a} and  
(using \cite{Graham95}) in
\cite{CoxGrahamMartin03}. 
The present paper is intended simply to provide an explicit
self-contained version (which can thus be put to
direct practical use in the study of further generalisations of the 
Temperley--Lieb algebra
{\em not} amenable to the methods of \cite{CoxGrahamMartin03}, 
that we shall report on elsewhere). 
The main results are Definition~\ref{pascal2blob} 
(a map from walks on the Pascal
triangle to `reduced' words in the blob algebra); 
Proposition~\ref{2} and Proposition~\ref{main bases}
(which together show that these words aggregate naturally into 
bases for certain blob algebra modules, including the regular
module). 


The blob algebra $b_n$ is usually defined in terms of a certain basis
of diagrams and their compositions (\cite{MartinSaleur94a}, and see
later), from which it derives its name.  
This blob algebra is isomorphic to an algebra defined by a presentation.  
There is a map which takes each generator in the presentation to a diagram,
which, it is fairly easy to see, provides a surjective algebra
homomorphism. 
The isomorphism follows from the vanishing of the kernel of this
homomorphism, which is somewhat harder to exhibit. 
A direct way to see it, however, 
is from a co-enumeration of the diagram basis
with a certain set of words, which is then shown to be spanning.
(The image of any spanning set under a surjective map is again
spanning, and hence, if of order no greater than the 
rank of the diagram algebra, independent. Then the original set is
independent, and hence a basis. We have a surjective homomorphism
between free modules of equal rank, and hence an isomorphism.)
This co-enumeration 
uses the natural tower structure $b_i \subseteq b_{i+1}$, and 
follows from the bra-ket formalism 
introduced in 
\cite{MartinSaleur94a}, via a generalisation of the corresponding 
proof of isomorphism between the diagram and presentation forms 
 \cite[\S6.4-6.5]{Martin91}
for the Temperley--Lieb algebra. 
In this way one gets bases not only for the algebra, but also for
Specht modules (cf. \cite{JamesKerber81}), 
and various double sided ideals, in a way which 
helps exhibit the structure of the algebra. 
Martin and Saleur \cite{MartinSaleur94a} 
eschewed details of this argument in the
interest of brevity, 
providing an explicit enumeration only on the diagram side. 
On the presentational side
a brief summary can be found in  \cite{CoxGrahamMartin03}
(using Graham \cite{Graham95} which,
motivated by the study of Hecke algebras,
extends the ideas of  \cite[\S6.5]{Martin91}
to an axiomatic framework broad enough to include
the blob algebra).
Here we provide an explicit 
self-contained
version on the presentation side, 
paralleling the use 
in \cite{MartinSaleur94a}
of an  `enriched' Pascal triangle. 

\ignore{{
In addition to proving isomorphism and 
constructing word bases for
a number of interesting blob algebra modules (see later)
\footnote{A more abstract approach, also extending the ideas in
  \cite[\S6.5]{Martin91}, but 
which does not provide this
  facility, is discussed in \cite{Graham95}
(see \cite{CoxGrahamMartin03} for a summary).}, 
the approach here also suggests a novel way to access  some
interesting combinatorial/representation theoretic results 
associated to the Pascal triangle \cite{James78}
(a point which will be pursued in detail in separate works 
\cite{MarshMartin05}). 
}}
\newcommand{\bb}{b}
\newcommand{\Uall}{U^e}
\subsection{The algebra $\bb_n$ by presentation}



\newcommand{\vargamma}{\gamma}%


For $K$ a ring, $x$ an invertible element in $K$, 
$q=x^2$, and $\vargamma,\vardelta \in K$, 
define $b_n^K$ to be the unital $K$--algebra with 
generator set $\Uall = \{ e,\U_1,\ldots,\U_{n-1} \}$ and relations
\begin{eqnarray} 
 \U_i \U_i &=& (q+q^{-1}) \U_i     \label{TL001} \\
 \U_i \U_{i\pm 1} \U_i &=& \U_i    \label{TL002} \\
 \U_i \U_j &=& \U_j \U_i \hspace{1in} \mbox{($|i-j|\neq 1$)}  \label{TL003}
\end{eqnarray}
\begin{eqnarray} 
 \U_1 e \U_1 &=&  \vargamma \U_1   \label{TL004} \\
 e e &=& \vardelta e \\
 \U_i e &=& e \U_i   \hspace{1in} \mbox{($i \neq 1$)} .   
\label{TL006}
\end{eqnarray}


It will be evident that $e$ can be rescaled to change $\vargamma$ and
$\vardelta$ by the same factor. 
Thus, if we require that $\vardelta$ is invertible, then we might as well
replace it by 1. 
(This brings us to the original two--parameter definition of the algebra.) 
We take $b_0^K = K\{1 \}$, 
and to be clear, note that $b_1^K = K\{1,e \}$. 


Note from the form of the relations that we have 
\prl(oppo) 
Let $\left( b_n^K \right)^{op}$ be the opposite algebra of $b_n^K$. 
Then there is an algebra isomorphism
\[  
b_n^K \cong \left( b_n^K \right)^{op}
\]
which fixes the generators. \Qed
\end{pr}


\noindent
For $k$ a field which is a $K$--algebra define $k$--algebra 
$b_n = k \otimes_K b_n^K$.





\prl(base, n=2) For $n=2$ 
\[
b_2^K = K \{ 1,e,\U_1, e\U_1, \U_1 e, e \U_1 e   \}
\]
(with the given set independent). 
\end{pr}
{\em Proof:} It will be evident that this set of words is spanning. 
(It is left as an exercise to show independence.) 
\Qed

\prl(finite bn)
Set   $\U_0 = e$. 
Then   
for $n \geq 1 $ 
\eql(bnbn-1)
b_{n}^K = b_{n-1}^K + b_{n-1}^K \U_{n-1} b_{n-1}^K . 
\eq
\end{pr}
{\em Proof:} By induction on $n$. 
The case $n=1$ of (\ref{bnbn-1}) is clear, 
and the case $n=2$ follows from proposition~\ref{base, n=2}. 
Now suppose true at level $n-1$, and consider $n>2$. 
Trivially 
\[
b_{n}^K = b_{n-1}^K + b_{n-1}^K \U_{n-1} b_{n-1}^K 
    + b_{n-1}^K \U_{n-1} b_{n-1}^K \U_{n-1} b_{n-1}^K + \ldots
\]
but by assumption
\eql(lem0)
\U_{n-1} b_{n-1}^K \U_{n-1} = 
\U_{n-1} b_{n-2}^K \U_{n-1} + \U_{n-1} b_{n-2}^K \U_{n-2} b_{n-2}^K \U_{n-1}
\stackrel{(\ref{TL002})}{=} \U_{n-1} b_{n-2}^K
\eq
\hfill \Qed

\prl(finbn2)
\eql(b2K)
\U_1 b_2^K \U_1 = ([2] K + \vargamma K )  \U_1 b_0^K 
\eq
and for $n \geq 3 $
\eql(bnK)
\U_{n-1} b_{n}^K \U_{n-1} = \U_{n-1} b_{n-2}^K
\eq
\end{pr}
{\em Proof:} (\ref{b2K}) follows from proposition~\ref{base, n=2} and
the defining relations. 
By proposition~\ref{finite bn}
\[
\U_{n-1} b_{n}^K \U_{n-1} = 
\U_{n-1} b_{n-1}^K \U_{n-1} + \U_{n-1} b_{n-1}^K \U_{n-1} b_{n-1}^K \U_{n-1}
\stackrel{(\ref{lem0})}{=} \U_{n-1} b_{n-2}^K
\]
\hfill \Qed


\medskip

In \cite[Proposition~2]{MartinSaleur94a} 
an explicit enumeration of the diagram basis of each blob
{\em diagram} algebra is given. 
Indeed this basis is put in explicit bijection with the set of pairs
of walks to the same location in level $n$ of Pascal's triangle. 
Here we do the same thing for bases for each algebra $b_n$ as defined
above.
We proceed by constructing bases for a series of ideals, the last of
which is $b_n$ itself.  

\ignore{

It is known \cite{MartinSaleur94a} that the representation theory of
$b_n$ falls into one of three distinct categories, 
depending on the number of integer values of $a$ for which 
$$ \vargamma [a]_q = \vardelta [a-1]_q $$
(where $[a]_q$ is the usual $q$--number). 
If there  is no solution then $b_n$ is semisimple. 
Accordingly it is convenient to reparameterize into one of the
following forms: 
\newline
$\vargamma=\frac{[\ym-1]}{[\ym]}$, 
$\vardelta=1$;
\newline
$\vargamma=[\ym -1]$, 
$\vardelta=[\ym]$;
\newline
$\vargamma=q^{\ym-1} - q^{-\ym+1}$, 
$\vardelta=q^{\ym} - q^{-\ym}$. 

The latter has the mild disadvantage that it is not a simple rescaling in
case $q-q^{-1} = 0$. 
Note, however, that provided $m$ is integer (which includes all the
interesting cases) this form has a lattice in $b_n^{\Z[q,q^{-1}]}$ 
(i.e. $\vargamma,\vardelta $  lie in $\Z[q,q^{-1}] $).

}




\section{Preliminaries on words in $\Uall$}
\subsection{Preliminary construction of ideals}
\newcommand{\Ue}{U^e}
\newcommand{\bdsi}{{\mathcal I}}

Let $\U=\{ U_1,U_2,\ldots,U_{n-1} \}$, so $\Ue=\{ e \} \cup \U$. 
For $S$ a set let $\Powerset(S)$ denote its power set. 
For $m\geq 0$, $m \equiv n$ (mod. 2) define
$$ \E{m}{(n)}   
= U_1 U_3 \ldots U_{n-m-1} 
\hspace{1in}
\Er{m}{(n)} = U_{n-1} U_{n-3} \ldots U_{m+1}  
$$
(take $\E{n}{(n)} =\Er{n}{(n)} = 1$). 
When $n$ is fixed we may simply write $\E{m}{}$ for $\E{m}{(n)}$, 
but note that
\eql(Eequiv) \E{m+1}{(n+1)} = \E{m}{(n)} . \eq
For $m > 0$, $m \equiv n$ (mod. 2) define
$$ \E{m+}{(n)} =\E{m+}{} 
   = \; \E{m}{(n)} \; e U_2 U_4 \ldots U_{n-m} \; \E{m}{(n)} . $$
\begin{de}
Define double sided ideals in $b_n$ by 
\[
\bdsi_m = \; b_n \; \E{m}{(n)} \; b_n  \hspace{1in}
\bdsi_{m}^+ = \; b_n \; \E{m+}{(n)} \; b_n  .  
\]
\end{de}
Thus $\bdsi_n = b_n$. 
%
It follows from proposition~\ref{finbn2} that
\prl(bss1)
For $m>0$
\[
\Er{m}{(n)} b_n \Er{m}{(n)} = \Er{m}{(n)} b_{m} 
\]
(for $m=0$ a similar result holds, but one must take care with the ring, 
as in proposition~\ref{finbn2} ---
at least one of $[2]$, $\gamma$ must be invertible).
\end{pr}
\begin{de}
Define
$\U^2 \subseteq {\Powerset}(\U)$ 
as the maximal subset 
such that $U_i,U_j \in V \in \U^2$
implies $i-j \not\in \{\pm 1 \}$ (i.e. the set of commutative subsets of $\U$). 
\end{de}
\prl(ideals)
\ignore{
Let $V,W,W' \in \U^2$ 
with $|V|>|W|=|W'|$. Then 
$$ b_n \left( \prod_{w\in W} w \right)  b_n =  
b_n \left( \prod_{w\in W'} w \right)  b_n $$
$$ b_n \left( \prod_{w\in V} w \right)  b_n \subset  
b_n \left( \prod_{w\in W} w \right)  b_n $$
}
Let $W \in \U^2$. Then 
$$ b_n \left( \prod_{w\in W} w \right)  b_n =  
\bdsi_{n-2|W|} $$
$$  \bdsi_{m} \subset   \bdsi_{m+2} $$
$$ \bdsi_{m}^+ \subset \bdsi_{m}
\hspace{1in}
\vargamma \bdsi_{m} \subset   \bdsi_{m+2}^+ $$
\end{pr}
The proofs are elementary applications of relations~(\ref{TL002})
(and, in the last case (\ref{TL004}), which provides the only means
to reduce out the `last' $e$ in any word). 

\newcommand{\bdelta}{\Delta}

Accordingly, for $i \equiv n$ (mod. 2), $i \geq 0$, define quotient algebras
\[
b_n^i = b_n / \bdsi_i
\hspace{1in}
b_n^{-i} = b_n / ( \bdsi_i^+ \cup \bdsi_{i-2} )
\]
Note that, for $n \geq 2$, $b_n^{n-2}$ has basis $\{ 1 , e \}$ 
(consider iterating proposition~\ref{finite bn} for example). 
Similarly, for $n-2 > 0$, 
$\Er{n-2}{(n)} b_n^{n-4} \Er{n-2}{(n)}$ 
has basis $\{ U_{n-1} , eU_{n-1} \}=\{ 1 , e \}\Er{n-2}{(n)}$ (as it were). 
Indeed, provided that $n-2r >0$,  
$\Er{n-2r}{(n)} b_n^{n-2r-2} \Er{n-2r}{(n)}$ 
has basis $\{ 1 , e \}\Er{n-2r}{(n)}$. 

For $n$ even, 
define left $b_n$-module $\bdelta_{0}(n) = b_n \E{0}{(n)} $.  
For $n$ odd, 
define left $b_n$-module $\bdelta_{-1}(n) = b_n \E{1+}{(n)} $,
and $\bdelta_{1}(n) = b_n \E{1}{(n)} $ mod. $\bdelta_{-1}(n)$.  
Define left $b_n$-module $\bdelta_{i+2}(n) $ 
as the restriction of the $b_n^{-i}$-module $b_n^{-i} \E{i+2}{(n)}$.
Define left $b_n$-module $\bdelta_{-(i+2)}(n) $ 
as the restriction of the $b_n^{i}$-module $b_n^{i} \E{(i+2)+}{(n)}$.

\subsection{Notations and identities}

Define
$$ U_{i\desc j} = U_i U_{i-1} \ldots U_j \hspace{1in} (i\geq j) $$
$$ U_{i\ddesc j} = U_i U_{i-2} \ldots U_j \hspace{1in} (i- j \in 2\N) $$
(and if the argument condition is violated we will take any such
product to evaluate to 1). 
NB the following elementary identities
\eql(Eident1) \E{0}(i)U_{i+1\desc 1} = \E{0}(i+2)  \eq
\eql(Eident1+) \E{1+}(i)U_{i+1\desc 1} = U_{i+1} \E{1+}(i+2)  \eq
\eql(chain1) U_{j\desc 1}U_{k\desc 1}
              = U_{j\desc 1}U_{k\desc 3} 
              = U_{j\desc 2}U_{k\desc 4} U_1 U_3
\hspace{1in} (k>j\geq 1) \eq
\eql(chain2) U_{j_1\desc 1}U_{j_2\desc 1}U_{j_3\desc 1}...
              = U_{j_1\desc 1}U_{j_2\desc 3} U_{j_3\desc 5}...
\hspace{1in} (j_3>j_2>j_1 \geq 1) \eq
$$
              = U_{j_1\desc 2}U_{j_2\desc 4} U_{j_3\desc 6}...U_{1}U_{3}U_{5}...
\hspace{1in} (j_i \geq 2i-1) $$
\eql(chain3) U_{2k\desc 1} U_{2j\ddesc 2} = U_{2j\ddesc 2}
\hspace{1in} (j \geq k) \eq
\eql(chain4) U_{2j-1\ddesc 1} U_{2k\desc 1} = U_{2j-1\ddesc 1}
\hspace{1in} (j \geq k) \eq

\subsection{Word reduction}

\newcommand{\ared}{algebra reduced}%

Note from (\ref{TL001}-\ref{TL006}) 
that every relation which shortens a word introduces a scalar
factor from $K$ (but that this factor may be 1). 
We call a word {\em \ared} if it cannot be expressed 
as a product of a scalar in $K\setminus\{1 \}$ times another word. 
(Thus the $K$-span of \ared\ words is the whole algebra
\cite{Bergman}.)
For example, $U_1 U_2 U_1$ is \ared. 

\prl(redux)

(a) Word $w \in b^K_{n-1}$ is \ared\ iff $wU_{n} \in b^K_{n+1}$ is
\ared.

(b) Word  $w \in b^K_{n-1}$ is \ared\ iff 
$wU_{n-1}U_{n-2}\ldots U_1 \in b^K_{n}$ is \ared.

(c) A word of form $wU_{n-2\ddesc 1} \in b^K_{n-1}$ is \ared\ iff 
$wU_{n-2\ddesc 1}eU_{n-1\ddesc 2}U_{n-2\ddesc 1} \in b^K_{n}$ is \ared.

\end{pr}
{\em Proof:}
The first claim follows from the commutation of 
$w$ with $U_{n}$ in $b^K_{n+1}$ 
(noting that any word which is \ared\ in $b^K_{n+1}$ but expressible in
$b^K_{n-1}$ is \ared\ in $b^K_{n-1}$). 
The second then follows since 
$(w U_{n})(U_{n-1}U_{n-2}\ldots U_1) U_2 U_3 \ldots U_{n} = w U_{n}$.
For the third note that a word of the form $wU_{n-2\ddesc 1}$ 
is \ared\ if and only if $wU_{n-2\ddesc 1}e$ is \ared.
This in turn is \ared\ if and only if 
$wU_{n-2\ddesc 1}eU_{n-1\ddesc 2}U_{n-2\ddesc 1}$ is \ared, since 
$U_n (wU_{n-2\ddesc 1}eU_{n-1\ddesc 2}U_{n-2\ddesc 1})
      U_{n-1\ddesc 2} U_{n\ddesc 3}
= wU_{n-2\ddesc 1}e U_n$, 
whereupon we can use (a) again 
(NB, the last identity is verified in $b^K_{n+1}$).
\Qed

\section{The Pascal triangle}
We now associate certain elements of $b_n$ to descending paths 
of length $n$ on the Pascal triangle. 
\subsection{Paths on the Pascal triangle}
Label vertices (positions) on the Pascal triangle by pairs of numbers
giving level (row) and
weight (column):
\[
\input{./xfig/PascalXY.eepic}
\]
Label edges by vertex pairs: $((n,m),(n+1,m\pm 1))$. 

For $m \in \Z\setminus\{ 0 \}$, $i \in \N$,  
define 
$$ m\out i = \left\{ \begin{array}{ll} m + i & \mbox{ in case $m>0$} \\
                                       m - i & \mbox{ otherwise.} 
\end{array} \right. $$ 
Define $m\back i = m \out -i$ similarly 
(except $m\back i =0$ in case $i>|m|$). 


\del(sn snm)
Let $S_n$ denote the set of walks from level 0 to level $n$ (any
weight) on the Pascal triangle; and $S_{n,m}$ the subset to weight $m$. 
\end{de}
It will be evident that  $|S_n|=2^n$.
For $p \in S_n$ we will write $p_i$ for the $i^{th}$ {\em edge} of $p$. 

There are various ways of specifying a particular $p \in S_n$. 
In particular, let $\sigma(p) =(\sigma(p)_0, \sigma(p)_1, \ldots)$ 
be the encoding of $p \in S_n$ as a sequence of weights. 
For example, $(0,1,0)$. 
The edges $p_i$ are then just the sequence of adjacent pairs from this
sequence. 
More robustly, we may specify a walk, or part of a walk, as a sequence
of (level,weight) pairs. (Our example becomes $((0,0),(1,1),(2,0))$.)
Then each edge is a pair of such pairs. 
 
\medskip
\subsection{Words and paths}
Associate words in the generators $\Ue$ (of $b_{\infty}$) 
to edges on the Pascal
triangle as:
\begin{eqnarray}
w(((n, |m|),(n+1, |m|+1)))&=&1  \label{p2w1} \\
w(((n,-|m|),(n+1,-|m|-1)))&=&1 \\
w(((n,m),(n+1,m \back 1)))&=& U_{n} U_{n-1} ... U_{1}  \label{p2w3} \\
w(((n,0),(n+1,+1)))&=&  
                        e U_{n\ddesc 2} U_{n-1\ddesc 1} \label{p2w4} \\
w(((n,0),(n+1,-1)))&=&1  \label{p2w5}
\end{eqnarray}

Each descending path $p$ on the Pascal triangle may be described 
by the sequence
$p_1,p_2,...$ of edges passed through. 
\del(pascal2blob)
For any $n$ define
\[
w: S_n \rightarrow \langle \Ue \rangle \; \in b^K_n
\]
as follows. 
The word $w(p)$ associated to path $p$ is the product
of words associated to this sequence of edges by (\ref{p2w1}--\ref{p2w5}), 
written from left to
right: $w(p)=w(p_1) w(p_2)...$. 
\end{de}
The first several such are given in table~\ref{PascalBase}. 
\begin{figure}
\input{./xfig/PascalBase.eepic}
\caption{\label{PascalBase} Words associated to edges, and hence
  paths, on the Pascal triangle. 
  Unlabelled edges have $w(\mbox{edge})=1$.}
\end{figure}

Note that if we write $w(p)=x \in b_n$ we mean the identity in $b_n$,
not necessarily as words. 

\del(snm1)
The set of words $w(p)$ associated to paths $p \in S_{n,m}$ 
(i.e. starting at vertex
$(0,0)$ and terminating at vertex $(n,m)$), is denoted $S_{(n,m)}$.  
\end{de}
\prl(1)   
\eql(p1e1)
S_{(n,m)} \subset b^K_n  \E{|m|}(n) .   
\eq

For $m>0$
\eql(p1e1+)
S_{(n,m)} \subset b^K_n \E{m+}(n)
\eq

For $m \leq 0$, $\E{m}(n)  \in S_{(n,m)}$; 
for $m   >  0$, $\E{m+}(n) \in S_{(n,m)}$.
\end{pr}
%
%
{\em Proof:} 
(Of (\ref{p1e1})): Suppose true at level $n-1$ 
(the base case is trivial). 
For $m=0$ we have 
$$
S_{(n,0)} = \{ xU_{n-1\desc 1} | x \in S_{(n-1, 1)} \} \cup 
            \{ xU_{n-1\desc 1} | x \in S_{(n-1,-1)} \}
$$
But $S_{(n-1,-1)} \subset b^K_{n-1} \E{1}(n-1)$ by the inductive
hypothesis, and 
$\E{1}(n-1) U_{n-1\desc 1} = 
 \E{0}(n-2) U_{n-1\desc 1} = \E{0}(n)$ 
so the second subset lies in $b^K_n \E{0}(n)$
(and similarly for the $S_{(n-1,1)}$  part).  
For $m \neq 0,1$ we have 
$$S_{(n,m)} = S_{(n-1,m\back 1)} \cup 
          \{ xU_{n-1\desc 1} | x \in S_{(n-1,m\out 1)} \}$$
The first subset obeys (\ref{p1e1}) by equation(\ref{Eequiv}) 
and the inductive hypothesis, 
the second by equation(\ref{Eident1}). 
For $m=1$ ((\ref{p1e1+})):
$$
S_{(n,1)} = \{ xU_{n-1\desc 1} \; | \; x \in S_{(n-1, 2)} \} \cup 
            \{ xeU_{n-1\ddesc 2}U_{n-2\ddesc 1} \; | \; x \in S_{(n-1, 0)} \}
$$
and we may proceed similarly using (\ref{Eident1}) and (\ref{Eident1+}).
For $m>1$ ((\ref{p1e1+})) use (\ref{Eident1+}) similarly.

(Of content claim): 
In case $m=0$ note that 
$w(0 \; -\! 1 \; 0 \; -\! 1\ldots 0) = \E{0}(n)$ by (\ref{Eident1}). 
In case $m=1$ note then that 
$w(0 \; -1 \; 0 \; -1\ldots 0 1) = \E{1+}(n)$ by (\ref{p2w4}). 
All the other cases follow by observing that every point on the Pascal
triangle can be reached by a path of the form 
$0 \; -\! 1 \; 0 \; -\! 1\ldots 0 \; -\! 1\; -\! 2\ldots \; -\! m \;$
or
$\; 0 \; -\! 1 \; 0 \; -\! 1\ldots 012\ldots m$, 
for which the last $m$ 
factors in $w(p)$ do not change the word. This $w(p)$ is thus the
required word (noting (\ref{Eequiv})). 
\Qed


\label{ss incl variant}

It follows that we may replace (\ref{p2w3}) by 
$$w(((n,m),(n+1,m \back 1)) = U_{n\desc n-|m\back 1|} $$
changing the word $w(p)$ only by an algebra equivalence. 
We will use the two forms interchangeably in what follows, 
unless the length of words
is an issue (in which case we will take the latter 
form, this being {\em never longer} than the former). 

It then follows from Proposition~\ref{redux} that
\prl(redux ad)
For all $n,m$: 
$S_{(n,m)}$ is a set of \ared\ words;
\\
 $S_{(n,m)}$ can be expressed in the form 
\[
S_{(n,m)} = \{ w \E{-m}(n) \; | \; w \in s(n,m) \}  \qquad \mbox{($m \leq 0$)}
\]
\[
S_{(n,m)} =  \{ w \E{m+}(n) \; | \; w \in s(n,m) \}  \qquad \mbox{($m > 0$)}
\]
where in each case $s(n,m)$ is some set of \ared\ words. 
\Qed
\end{pr}

\medskip

\section{Ideals, modules and bases}

Consider the action of $U_{i}$ ($i \in \{1,2,..,n-1 \}$) 
from the left on each element $w(p)$ of $S_{(n,m)}$. 
For this purpose the 
most significant part of the {\em path} $p$ 
is the neighbourhood of the $i^{th}$ vertex, 
as we will see. 
This part may be specified by expanding $p$ in the form 
$$p=p'(i-1,l)(i,m)(i+1,n)p'' = p'(l,m,n)p'' $$ 
for some weight triple $(l,m,n)$. 
(This is an abuse of notation, since weight $(i-1,l)$ is also a part of
the final edge in subpath $p'$, but it is still useful.) 
In this sense we may write $w(p)=w(p')w((l,m,n))w(p'')$  
with $w(p') \in b_{i-1}$ (and hence commuting with $U_i$). 
For example 
$$ w(p'(l,l\back 1, l)p'') = \; w(p') \; U_{i-1\desc 1} \; w(p'') $$ 
and hence
\prl(diamond1) For $|l|>1$ 
$$ U_i \; w(p'(l,l\back 1, l)p'') \; = w(p'(l,l\out 1,l)p'') $$
\Qed
\end{pr}

It will be convenient to picture the difference between the two walks
here as follows: 
\[
\includegraphics{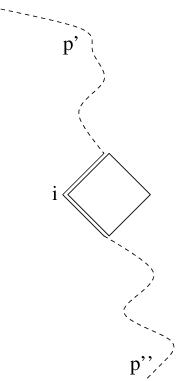}
\]

Let us generalise the notation of the 
triple $(l,m,n)$ (whose middle element is the
weight at the $i^{th}$ vertex) to any sequence of weights of
consecutive vertices (NB, the position of the $i^{th}$ vertex
must be indicated in some way). For example
\prl(diamond0) For $l \in \N$ 
$$ U_i \; 
w(p'(0(-\! 1 \; -\! 2)^l \; -\! 1 \; \overbrace{0}^{i^{th}} \; 1)p'') 
\; =  w(p'(0(1 \; 2)^l 1 \; 2 \; 1)p'') $$
\end{pr}
The picture here is:
\includegraphics{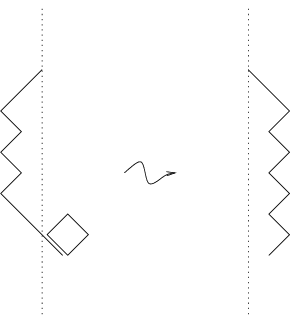}

\noindent
{\em Proof:} Note that $i$ is even. Put $j=i-2l-2$, then  
simply applying (\ref{p2w1}-\ref{p2w5}) we get
$$  w(p'(0(-\! 1 \; -\! 2)^l \; -\! 1 \; 0 \; 1)p'')  \hspace{9.8cm}$$
$$ =  w(p') \; 
U_{j\desc 1} U_{j+2\desc 1} \ldots U_{j+2l\desc 1}  U_{j+2l+1\desc 1} 
e U_{2}U_{4}\ldots U_{j+2l+2} U_{1} U_{3} \ldots  U_{j+2l+1} \; w(p'')
. \hspace{.1cm} $$ 
Thus 
$$ U_i \; w(p'(0(-\! 1 \; -\! 2)^l \; -\! 1 \; 0 \; 1)p'')  \hspace{8cm}$$
$$ = U_i \; w(p') \; 
U_{j\desc 1} U_{j+2\desc 1} \ldots U_{j+2l\desc 1}  U_{j+2l+1\desc 1} 
e U_{2}U_{4}\ldots U_{j+2l+2} U_{1} U_{3} \ldots  U_{j+2l+1} \; w(p'')
\hspace{1cm} $$ 
$$ = w(p') U_{j\desc 1} U_{j+2\desc 1} \ldots U_{j+2l\desc 1}
U_{j+2l+2\desc 1} 
e U_{2}U_{4}\ldots  U_{j+2l+2} U_{1} U_{3} \ldots  U_{j+2l+1} \; w(p'') . $$ 
But using  equation(\ref{chain3}) repeatedly
$$  U_{j\desc 1} U_{j+2\desc 1} \ldots U_{j+2l\desc 1}
U_{j+2l+2\desc 1} 
e U_{2}U_{4}\ldots  U_{j+2l+2} U_{1} U_{3} \ldots  U_{j+2l+1} 
$$ $$= e U_{2}U_{4}\ldots  U_{j+2l+2} U_{1} U_{3} \ldots  U_{j+2l+1} $$ 
while
$$ w(p'(0(1 \; 2)^l 1 \; 2 \; 1)p'') = 
\hspace{10cm} $$ $$
w(p') ( e U_2 U_4 \ldots U_{j-2} U_1 U_3 \ldots  U_{j-3}) 
U_{j\desc 1} U_{j+2\desc 1} \ldots U_{j+2l\desc 1} U_{j+2l+2\desc 1}
w(p'') $$
and the result follows from repeated application of equation(\ref{chain4}). 
\Qed


\prl(diamond012) For $l \in \N$ 
$$ U_i \; w(p'(0\; \overbrace{1}^{i^{th}}( 2 3)^l \; 2  \; 1)p'') 
\; =  w(p'(0 \; -\!\! 1 (0 \;  -\! 1)^l 0 \; 1)p'') $$
\end{pr}
The picture here is:
\[
\includegraphics{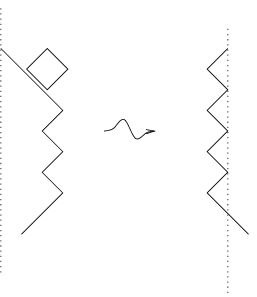}
\]
\noindent
{\em Proof:} 
The $l=0$ case follows from
proposition~\ref{diamond0} and the others by a simple iteration: 

If $p=p'(01232...)$
$$
U_{i}w(p) = 
    w(p')   (U_i) (eU_2\ldots U_{i-1} U_1 \ldots U_{i-2})
    U_{i+2}U_{i+1} U_{i}U_{i-1}..U_{1} w(p''')
$$
$$= 
    (U_{i+2}) 
    w(p')   (U_i) (eU_2\ldots U_{i-1} U_1 \ldots U_{i-2})
    U_{i+1} U_{i}  
    w(p''')$$ 
\eql(d012)= 
    (U_{i+2}) 
    w(p')   (U_i) (eU_2\ldots U_{i-1} U_{i+1}U_1 \ldots U_{i-2} U_{i})
    w(p''')  
=   (U_{i+2}) 
    w(p' (0 \; -\!1 \; 0 12 ...))  
\eq
and so on. 
\Qed


\prl(2)
If $m=0,1$ 
then $S_{(n,m)}$ spans $b_n S_{(n,m)}$;   

if $m=-1 $ then 
$S_{(n,m)}$ spans a  left subideal 
mod. $b_n S_{(n,1)} $. 
 
if $m\geq 2$ then 
$S_{(n,m)}$ spans a  left subideal mod. 
$b_n \E{m\back 2}(n) b_n $;

and if $m< -1$ then 
$S_{(n,m)}$ spans a  left subideal 
mod. $b_n S_{(n,-m)} \cup b_n \E{-(m\back 2)}(n) b_n$. 
\end{pr}
{\em Proof:} 
It is relatively straightforward to check that the action of $e$ on
$S_{(n,m)}$ stays within the indicated span.
Thus we concentrate on the action of $U_i$. 
Again consider the action of $U_i$ on $w(p)$, and characterise the
path in the neighbourhood of $i$ by the triple $(l,m,n)$. 

First suppose that one or both of the edges touching the ${i}^{th}$ vertex 
touches $m=0$. There are various cases:
\begin{itemize}\item
$(0,\pm 1,0)$: 
$$p=p'(i-1,0)(i,\pm 1)(i+1,0)p''$$ (NB, $i$ odd) 
then $w(p)=w(p')w((i-1,0)(i,\pm 1))U_{i}...U_{1}w(p'')$ 
so $$U_{i}w(p)=[2] w(p) . $$ 
\item
$(0, 1,  2)$: 
If $p''$ never turns over then $w(p'')=1$ and 
$U_{i}w(p) \in b_n \E{\pm(m\back 2)}(n)$ by
proposition~\ref{ideals}. 
If $p''$ never touches $m=1$ again then 
we may use equation(\ref{d012}) to equate to a case in which it never
turns over, and again  
$U_{i}w(p) \in b_n \E{\pm(m\back 2)}(n)$ by
proposition~\ref{ideals}. 
Otherwise we may use proposition~\ref{diamond012} (in combination with
proposition~\ref{diamond1}). 

Cases $(0,- 1, - 2)$, $(2,1,0)$ and $(-2,-1,0)$ are simpler, 
but with a similar strategy. 

\ignore{If $p''$ turns over at $j>i$ then
$w(p'')=U_{j}U_{j-1}..U_{i+1}U_{i}..U_{1} w(p''')$ so 
$U_{i}w(p) = 
    w(p')   (U_i) (eU_2\ldots U_{i-1} U_1 \ldots U_{i-2})
    U_{j}U_{j-1}...U_{i+1} U_{i}U_{i-1}..U_{1} w(p''')$ 
$= 
    (U_{j}U_{j-1}...U_{i+2}) 
    w(p')   (U_i) (eU_2\ldots U_{i-1} U_1 \ldots U_{i-2})
    U_{i+1} U_{i}  
    w(p''')$ 
$= 
    (U_{j}U_{j-1}...U_{i+2}) 
    w(p')   (U_i) (eU_2\ldots U_{i-1} U_{i+1}U_1 \ldots U_{i-2} U_{i})
    w(p''') . $
Either this is in  $b_n \E{\pm(m\back 2)}(n)$ by
proposition~\ref{ideals} (and we may stop) or 
$=   (U_{j}U_{j-1}...U_{i+2})  w(p')   (010 -1 -2 ...)  w(p''')$ 
where $p'''$ again turns over (at $j'>j$), whereupon we may repeat:
$=   (U_{j}U_{j-1}...U_{i+3})(U_{j'}U_{j'-1}...U_{i+4})  w(p')   (01010 -1 -2 ...)
    w(p'''')$ 
and so on. Eventually, if the process terminates no sooner, we reach a
path $(010101...010 -1 )$ (with a string of $U_l$s acting).  
}

\item
In case $(+1,0,-1)$ ($i$ even) then $w(p)=w(p')U_{i-1\desc 1} w(p'')$
so $$U_i w(p) = w(p')w((+1,+2,+1))w(t(p''))$$ where $t$ reflects that
part of its
argument up to the first $m=0$ in $m=0$. If $p$ ends at this first
zero or beyond we are
done. If $p$ ends before touching zero again then we changed the sign
of the end point of the walk (to positive) and we are done by the quotient. 

In case $(-1,0,+1)$ ($i$ even) every case is covered by some
combination of propositions \ref{diamond1} and \ref{diamond0}. 
\end{itemize}
If neither of the edges touching the ${i}^{th}$ vertex 
touches $m=0$, then:
\begin{itemize}
\item
If 
$$p=p'(i-1,m)(i,m\out 1)(i+1,m)p''$$ 
then $w(p)=w(p')U_{i}...U_{1}w(p'')$ 
so $$U_{i}w(p)=[2] w(p) . $$ 
\newline
If $p=p'(i-1,m)(i,m\back 1)(i+1,m)p''$ 
then $w(p)=w(p')U_{i-1}...U_{1}w(p'')$ 
so $$U_{i}w(p)= w(p'(i-1,m)(i,m\out 1)(i+1,m)p'') . $$ 
%
\newline
If $p=p'(i-1,m\back 2)(i,m\back 1)(i+1,m)p''$ 
then $w(p)=w(p')w(p'')$. 
If $p''$ never turns over then $w(p'')=1$ and  
$U_{i}w(p) \in b_n \prod_{j}U_j$ by
proposition~\ref{ideals}. 
If $p''$ turns over at $j>i$ then
$w(p'')=U_{j}U_{j-1}..U_{i+1}U_{i}..U_{1} w(p''')$ so 
$U_{i}w(p) = 
    ( U_{j}U_{j-1}...U_{i+2} ) w(p') U_{i}U_{i-1}..U_{1} w(p''')$ 
and 
$$ w(p')  U_{i}U_{i-1}..U_{1} w(p''') 
   = w(p' (i-1,m\back 2)(i,m\back 1)(i+1,m\back 2)...) . $$ 

%
We may partially order the set of walks which are identical except in
some interval where neither crosses $m=0$ 
(and hence the sign of $m$ never changes) by $p\geq p'$ if $|m|\geq |m'|$
throughout this interval. 
We have converted the action of $U_i$ on $w(p)$ 
in the case above to an action on 
a lower walk. We will return to this case shortly.
%

If $p=p'(i-1,m\out 2)(i,m\out 1)(i+1,m)p''$ then 
$w(p)=w(p')U_{i-1}...U_{1}U_{i}...U_{1} w(p'')$ 
so 
$U_{i}w(p) $ 
may  be reduced to the action of a string of $U_j$s (depending on
$p'$) on a lower walk, in a manner analogous to the case above. 


We have converted the action of $U_i$ on $w(p)$ 
in each of the two cases above to an action on 
a lower walk, thus we may apply an induction with one of the $m$
touching zero cases as base. 
\end{itemize}
\Qed

\del(dsia)
For $a \in b_n$ let $a^o$ denote the image under the opposite isomorphism.
Let $S^{2}_{(n,m)}$ denote the set of words 
$\{ a (\E{m+}(n)) b^o \; | \; a,b \in s{(n,m)}  \}$ 
(here if $m \leq 0$ then $m+$ means $-m$).
\end{de}
For example, reading from figure~\ref{PascalBase} we have 
\[
S_{(3,1)} = \{ U_1 e U_2 U_1, e (U_1 e U_2 U_1), U_2 (U_1 e U_2 U_1) \}
\]
so
\ignore{{
\[
S^{2u}_{(3,1)} =  \begin{array}{ccc} \{ \; 
U_1 e U_2 U_1  U_1 U_2 e U_1, & e U_1 e U_2 U_1  U_1 U_2 e U_1, & 
e U_2 U_1  U_1 U_2 e U_1, \\ \;
U_1 e U_2 U_1  U_1 U_2 e U_1 e, & e U_1 e U_2 U_1  U_1 U_2 e U_1 e, & 
e U_2 U_1  U_1 U_2 e U_1 e, \\
U_1 e U_2 U_1  U_1 U_2 e,  & e U_1 e U_2 U_1  U_1 U_2 e,  & 
e U_2 U_1  U_1 U_2 e  
\} \end{array} 
\]

The words in $S^{2u}_{(n,m)}$ are not reduced, 
but they all differ from reduced words
by the same scalar. Define  $S^{2}_{(n,m)}$ as the underlying set of 
reduced words. 

For example,
}}
\[
S^{2}_{(3,1)} =  \begin{array}{ccc} \{ \; 
U_1 e U_2  U_1, & e U_1 e U_2  U_1, & 
e U_2  U_1, \\
U_1 e U_2  U_1 e, & e U_1 e U_2  U_1 e, & 
e U_2  U_1 e, \\
U_1 e U_2 ,  & e U_1 e U_2 ,  & 
e U_2   
\} \end{array} 
\]

By Proposition~\ref{1} 
every $S_{(n,m)}$ contains an element invariant under the opposite
isomorphism 
($U_1 e U_2 U_1$ in our example)
so, noting that the above argument works analogously for
right ideals, 
from the definition and the last proposition we have
\prl(dsi1)
$S^{2}_{(n,m)}$ spans the double sided ideal it generates,
modulo the double sided ideals generated by all $S^2_{(n,m')}$ with
$|m'|<|m|$ and $m'=-m$ if $m<0$. 
\Qed
\end{pr}
Since $S^2_{(n,-n)} = \{ 1 \}$ we have that 
\eql(regular)
S^2_{(n)} :=  \bigcup_m S^2_{(n,m)}
\eq
spans $b_n$.

Finally, in the next section, we show linear independence. 
 
\section{Modules, bases and diagrams}

A blob diagram is a Temperley--Lieb diagram in which any line which
may be deformed isotopically to touch the western edge of the frame
may be decorated with a `blob'. The set of such diagrams with $n$
vertices on each of the northern and southern edge is denoted $B_n$. 
Two diagrams $d_1,d_2 \in B_n$ are `concatenated' by a juxtaposition
which identifies each southern vertex of $d_1$ with a northern
vertex of $d_2$ in the natural way.  
The blob ({\em diagram}) algebra,
here denoted  $b_n'$, 
 is defined in \cite{MartinSaleur94a}. 
In short, 
$b_n'$ has basis $B_n$, with  composition on this basis defined
as follows. The concatenation of $a,b \in B_n$ gives another diagram, 
except that this may contain some extra features. 
We interpret this pseudodiagram in $b_n'$ as follows. 
Firstly each internalised vertex is ignored. 
Each undecorated loop is removed and interpretted as a scalar factor
$\delta$; 
each occurence of a second blob on a given line is removed
and  interpretted (in our implementation) as a scalar factor $\delta_e$; 
each decorated loop is then removed and interpretted as a scalar factor
$\gamma$. (Note that after these removals the diagram will again lie
in $B_n$.)

The map  $\phi:b_n \rightarrow b_n'$ is given by 
\[ U_i \mapsto 
\raisebox{-.21in}{
\includegraphics{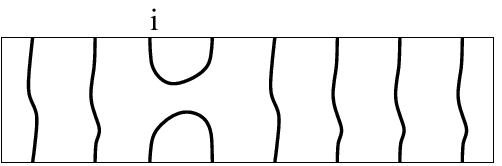}
}
\]
\[ e \mapsto 
\raisebox{-.21in}{
\includegraphics{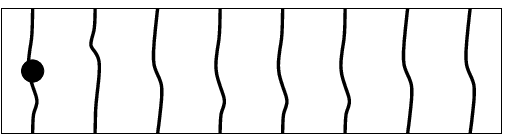}
}
\]
\prl(surj)
The map  $\phi:b_n \rightarrow b_n'$ is a surjective algebra homomorphism.
\end{pr}
{\em Proof:}
That this is an 
algebra homomorphism follows from a straightforward check of the
relations. For surjectivity, 
compare the image of $U_{i\desc 1}$ with the usual rule for
constructing half--diagrams using the Pascal triangle
\cite{MartinSaleur94a}. 
(Alternatively, just note that the images of the generators generate
$b_n'$.) 
\Qed

\medskip



Since the degree of $S^2_{(n)}$ conincides with the rank of $b_n'$,
the surjectivity of the map $b_n \rightarrow b_n'$
implies 
\prl(main)
For any $n$, $\phi$ defines an isomorphism 
$b_n \cong b_n'$.
\Qed
\end{pr}
And hence
\prl(main bases)
Every spanning set constructed in 
propositions~\ref{2} and \ref{dsi1}, and (\ref{regular}), is a basis
for the corresponding module. 
\Qed
\end{pr}

\[ \]
{\bf Acknowledgement.}
I thank 
R J Marsh for useful comments, and RJM and 
A~E~Parker for encouraging me to make these notes available.

\appendix
\section{\hspace{-20pt}ppendix}
\subsection{A compendium of related combinatorial facts}

Consider the `diamond' grid of side length $n$ (that is, the square grid
with $n+1 \times n+1$ vertices, oriented at $45^o$). 
Let $T_n$ denote the set of right-stepping walks from the set of vertices on
the centre vertical of this grid to the rightmost vertex (in bijection
with the set of walks from leftmost to rightmost which are symmetric
about the centre vertical). 

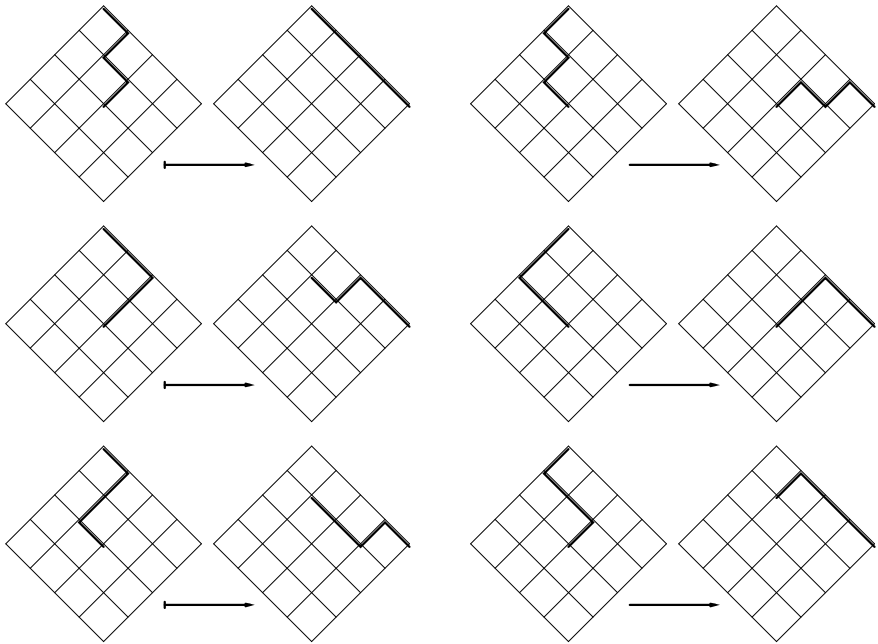
\begin{figure}
\input{./xfig/bigdiamond02.eepic}
\caption{\label{bigdiamond02} $\pi:S_{4,0} \hookrightarrow T_4$.}
\end{figure}

There is an obvious bijection between $T_n$ and $S_n$ got by rotating
through $90^o$. We wish to construct a {\em different bijection}. 
\ignore{
Starting from the endpoint of $p \in S_n$ parse the steps to obtain
steps in the image of $p$ also from last to first, as follows. 
Each step back ($m \rightarrow m\back 1$) parses to a step down; and
each step out parses to a step up {\em unless} the step out touches
$m=0$, in which case is parses to a step down. 
}
For $p \in S_n$, 
parse the edge sequence $\sigma(p)$ 
to a sequence $\pi(p)$ of elements from the set
$\{ N, S \}$ as follows: reading $\sigma(p)$ from left to right, 
\newline
if $|\sigma(p)_i | < |\sigma(p)_{i-1} | $ then  $\pi(p)_i=S$; 
\newline
if $(\sigma(p)_{i-1},\sigma(p)_{i}) \neq (0,1)$ and 
$|\sigma(p)_i | > |\sigma(p)_{i-1} | $ then  $\pi(p)_i=N$; 
\newline
if $(\sigma(p)_{i-1},\sigma(p)_{i}) = (0,1)$ then  $\pi(p)_i=S$. 
\newline
For example $\pi((0,1,0)) = (S,S)$. 
Each sequence $\pi(p)$ encodes a walk in $T_n$ by regarding $N$ as a
northeast step and $S$ as a southeast step, with the starting point
determined by the requirement that the finishing point is fixed for
all walks. 
Some examples are shown in figures~\ref{bigdiamond02} and~\ref{bigdiamond04}.
\begin{figure}
\input{./xfig/bigdiamond04.eepic}
\caption{\label{bigdiamond04} $\pi:S_{4,\pm2} \hookrightarrow T_4$.}
\end{figure}


Let `heights' on the diamond grid be measured from the lowest vertex
(height 0). Define a poset $(T_n,\geq)$ by $t\geq t'$ if,
reading from left to right, at each point the height of $t$ is $\geq$
the height of $t'$. 
(We will induce a poset $(S_n,\geq)$ from this, via $\pi^{-1}$.) 
\prl(bdw) 
$\pi$ is a bijection. 
$\pi(S_{n,m})$ is the subset of $T_n$ containing those walks whose
lowest point is at height $n-|m|+\frac{m+|m|}{2m}$. 
\end{pr}
The proof is elementary. 


Define a map $W:T_n \rightarrow b_n$ as follows.
Draw the walk $t$ together with
the maximal walk with lowest point $n-|m|$ ($m$ taken from $t=\pi(p)$ as
above). 
In each small diamond in the envelope so created write the
generator corresponding to that position (writing $e$s in the centre
`half--diamonds'). Now read off $W(t)$ from this picture from left to
right, top to bottom. 
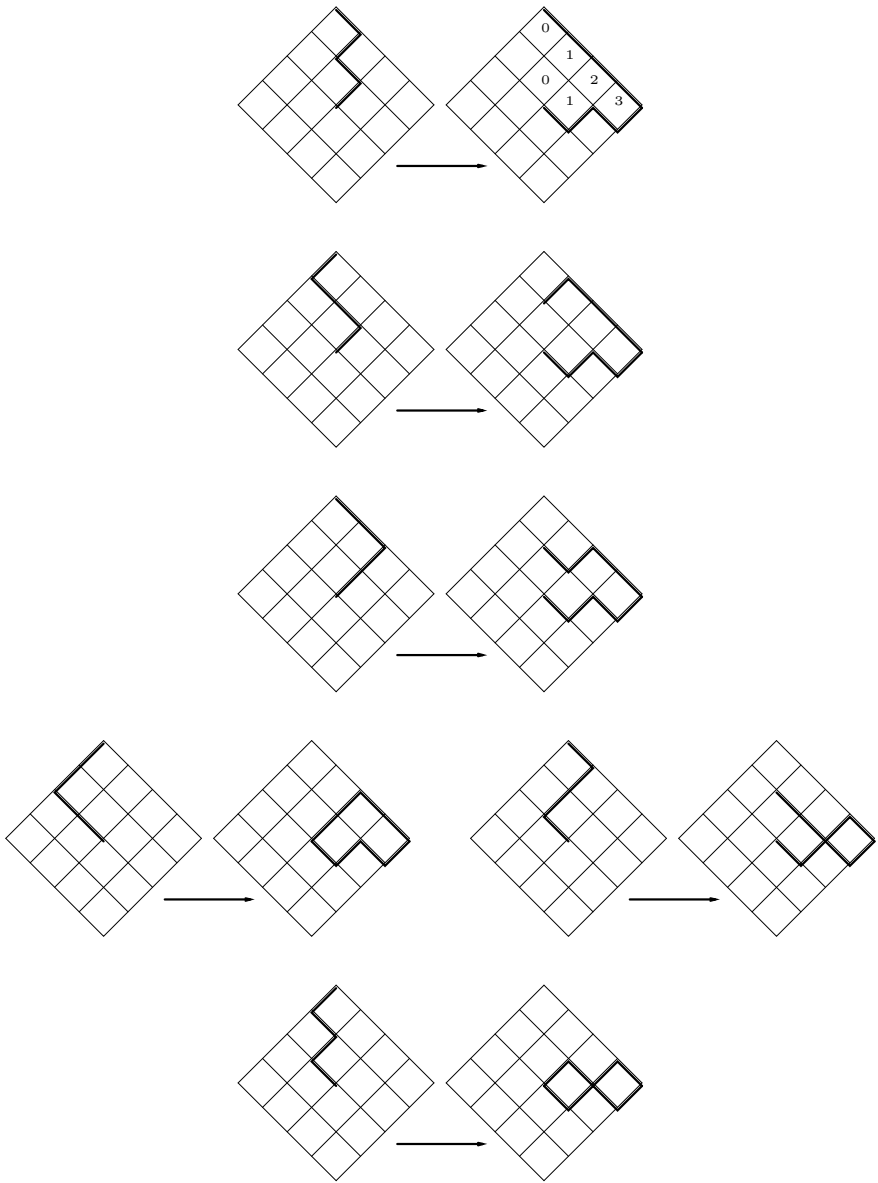
\begin{figure}
\input{./xfig/bigdiamond03.eepic}
\caption{\label{bigdiamond03} $\pi(S_{4,0})$ arranged in poset order. For
  the topmost case $w(p)=W(\pi(p))=e U_1 e U_2 U_1 U_3$. } 
\end{figure}
\begin{figure}
\input{./xfig/bigdiamond05.eepic}
\caption{\label{bigdiamond05} $\pi(S_{6,0})$ arranged in poset order. For
  the topmost case 
  $w(p)=W(\pi(p))=e U_1 e U_2 U_1 U_3 e U_2 U_4 U_1 U_3 U_5$. } 
\end{figure}
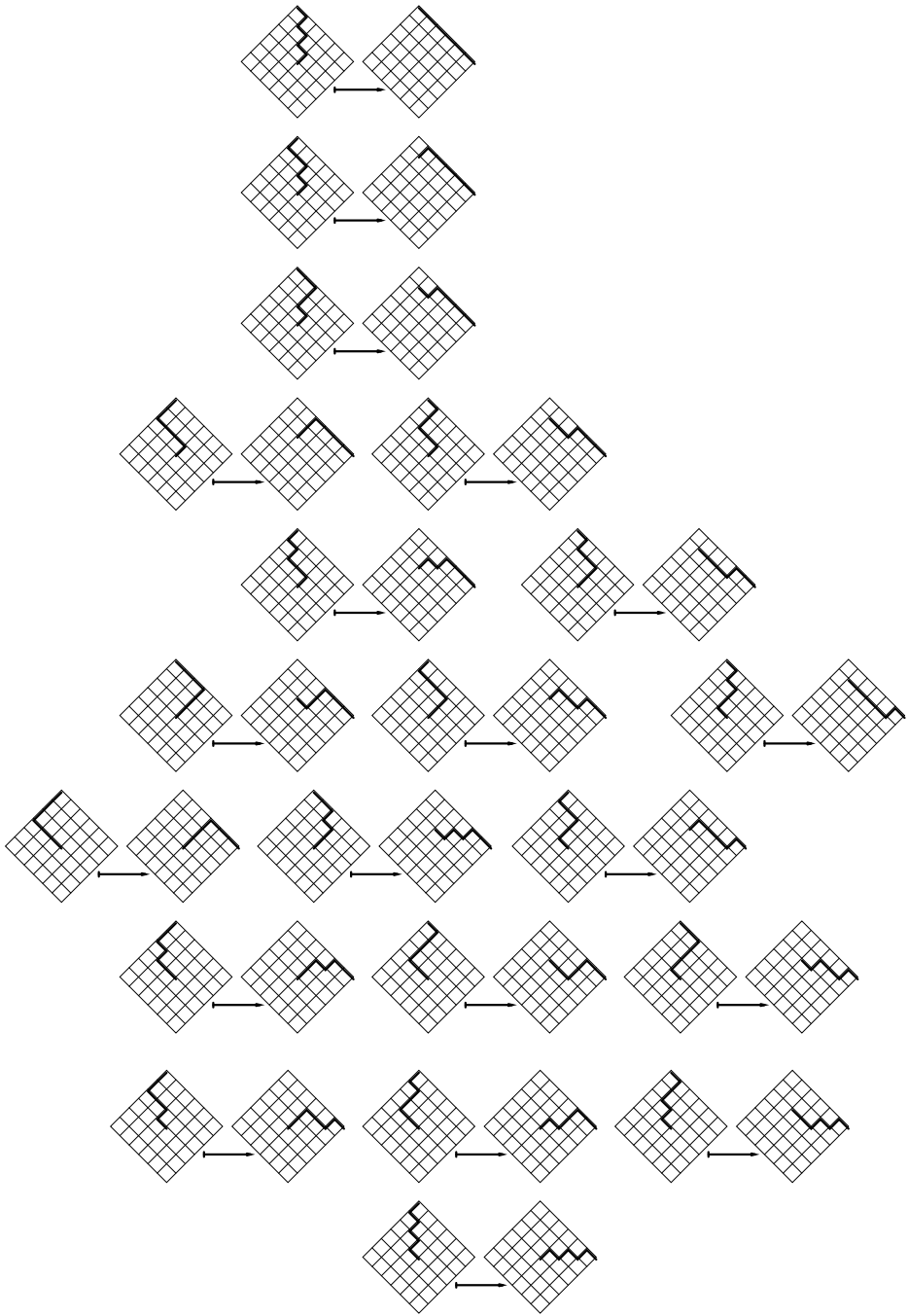
The walks from the example in figure~\ref{bigdiamond02} are shown
arranged in the partial order (top to bottom) in
figure~\ref{bigdiamond03}. The lower envelopes and indices for the
relevant $U_i$s have also been drawn (with $U_0=e$). 
(A bigger example, with $n=6$, is given in figure~\ref{bigdiamond05}.)
\prl(Ww)
$W(\pi(p))=w(p)$
\end{pr}
The proof is elementary.  





\subsection{Word basis: variant form}
\begin{figure}
\input{./xfig/PascalBase02.eepic}
\caption{\label{PascalBase02} Words associated to edges, and hence
  paths, on the Pascal triangle. 
  Unlabelled edges have $w(\mbox{edge})=1$.}
\end{figure}
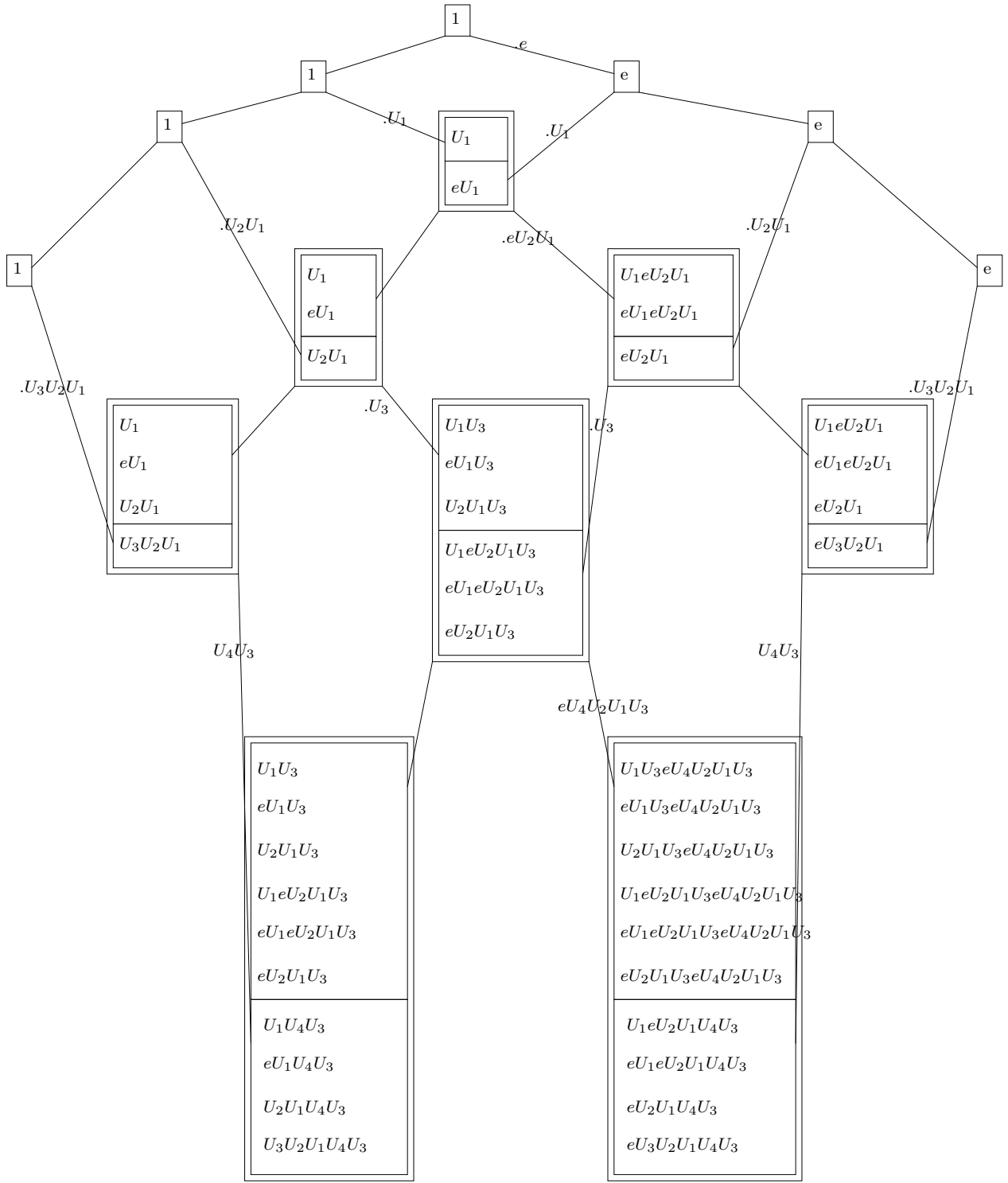

The variant form of the word set on the Pascal triangle 
(from the end of section~\ref{ss incl variant})
is shown in
figure~\ref{PascalBase02}. 



\bibliographystyle{amsplain}
\bibliography{new31,main,emma}

\end{document}

%% file: martinew.tex
\providecommand{\noglossaryignore}[1]{}
\newcommand{\globalglossaryentry}[3]{\makebox[1.5in][l]{\tt $\backslash${#1}} 
\makebox[1.1in][l]{{$#2$}} \makebox[2.5in][l]{{#3}}\newline} 
\newcommand{\newcommandabbreviation}[3]{\newcommand{#1}{#2}%
\noglossaryignore{\globalglossaryentry{#3}{#2}{}}}
\newcommand{\renewcommandabbreviation}[3]{\renewcommand{#1}{#2}%
\noglossaryignore{\globalglossaryentry{#3}{#2}{}}}
\newcommand{\newcommandmacro}[4]{\newcommand{#1}{#2}%
\noglossaryignore{\globalglossaryentry{#3}{#2}{#4}}}
\newcommand{\gge}[3]{\noglossaryignore{\globalglossaryentry{#1}{#2}{#3}}}
\newcommand{\myaddress}%
{\parbox{3in}{\footnotesize \begin{center} 
Mathematics Department, City University, \\  
Northampton Square, London EC1V 0HB, UK.\end{center}}}
    
\newcounter{minidef}[section]

\newcounter{minicapt}

\newtheorem{de}{Definition}     \newtheorem{pr}{Proposition}

\noglossaryignore{GREEK ETC.\newline}
\newcommandabbreviation{\e}{\epsilon}{e}        
\newcommandabbreviation{\lam}{\lambda}{lam}  
\newcommandabbreviation{\la}{\langle}{la}        
\newcommandabbreviation{\ran}{\rangle}{ran}
\newcommandabbreviation{\ha}{\#}{ha}             
\newcommandabbreviation{\rmap}{\rightarrow}{rmap}
\newcommandabbreviation{\aaa}{\alpha}{aaa}        
\newcommandabbreviation{\ab}{\alpha,\beta}{ab}
\newcommandabbreviation{\aab}{a(\ab )}{aab}       
\noglossaryignore{\newline RINGS\newline}
\newcommandabbreviation{\HH}{H \!\!\! I}{HH}               
\newcommandabbreviation{\C}{\mathbb C}{C}
\newcommandabbreviation{\N}{\mathbb N}{N}   
\newcommandabbreviation{\Z}{\mathbb Z}{Z}      
\renewcommandabbreviation{\Re}{\mathbb R}{Re}
\newcommandabbreviation{\R}{{\mathbb R}}{R}
\newcommandabbreviation{\Q}{\mathbb Q }{Q}
\renewcommandabbreviation{\H}{\mathbb H }{H}
\noglossaryignore{\newline SYMMETRIC GROUP\newline}
\def\Sym(#1){\Sigma(#1)}                   
\gge{Sym(-)}{\Sym(-)}{symmetric group on - objects}
\def\Sy(#1){\Sigma_{#1}}                   
\gge{Sy(-)}{\Sy(-)}{symmetric group irreducible -}
\def\sym(#1){\mbox{\LARGE s}(#1)}        
\gge{sym(-)}{\sym(-)}{symmetric group on - objects (variant)}
\def\sy(#1){\mbox{\LARGE s}({#1})}        
\gge{sy(-)}{\sy(-)}{symmetric group irreducible - (variant)}
\newcommandmacro{\cs}{\C \, \sy(n)}{cs}{symmetric group algebra over $\C$}
\noglossaryignore{\newline PARTITIONS/SETS\newline}
\newcommand{\Nset}[1]{\underline{#1}}
\gge{Nset\{-\}}{\Nset{-}}{set of natural numbers to -} 
\def\nset(#1){ \{ #1 \}_{ \underline{n} }} 
\gge{nset(-)}{\nset(-)}{a set $-\times\Nset$}
\def\ul(#1){_{\underline{#1}}}             
\gge{ul(-)}{{}\ul(-)}{subscript underline -}
\def\Ee(#1){{\bf E}_{#1}}                  
\gge{Ee(-)}{\Ee(-)}{set of equivalence relations on set -}
\def\Eee(#1){{\bf E}_{\{ #1 \}_{\underline{n}}}}   
\gge{Eee(-)}{\Eee(-)}{ditto for nset}
\def\Een(#1,#2){{\bf E}_{\{ #1 \}_{\underline{#2}}}}   
\def\Ssn(#1,#2){{\bf S}_{\{ #1 \}_{\underline{#2}}}}   
\def\Ss(#1){{\bf S}_{#1}}                  
\def\Sss(#1){{\bf S}_{\{ #1 \}_{\underline{n}}}}   
\def\bbc(#1){((\beta_1)(\beta_2)...(\beta_{#1}))}      
\newcommandmacro{\Ln}{{\Gamma}^{n}}{Ln}{large index set}
\newcommandmacro{\LnQ}{{\Gamma}^{n}_Q}{LnQ}{index set}
\newcommandmacro{\Zz}{\zeta}{Zz}{`shape' function}
\noglossaryignore{\newline PARTITION ALGEBRA\newline}
\def\ka(#1){\kappa_{#1}}                   
\def\Sm(#1){\Sigma_{#1}}                   
\newcommandmacro{\com}{\bullet}{com}{bullet composition}
\newcommandmacro{\enm}{\; e^n(\! m\! ) \;}{enm}{product of idempotents}
\def\Ai(#1){ A^{ #1 \cdot } }              
\def\Aij(#1,#2){ A^{ #1  #2 } }            
\newcommandmacro{\One}{\mbox{\bf $1 \!\!\! 1$}}{One}{algebra unit 1}
\newcommandmacro{\Bp}{B_p}{Bp}{partition basis}
\def\Bb(#1){B_p[#1]}                       
\def\Pp(#1){P_n[#1]}                       
\def\Ps(#1){P_n[#1] \! /}                  
\newcommandmacro{\Ph}{\hat{P}}{Ph}{P hat  algebra}
\def\Is(#1){\sim^{#1}}                     
\noglossaryignore{\newline MODULES\newline}
\def\Wm(#1){{\cal S}_{#1}}                 
\gge{Wm(-)}{\Wm(-)}{Weyl module with index -}
\def\wm(#1,#2){{}_{#1}{\cal S}_{#2}}       
\gge{wm(-1,-)}{\wm(-1,-)}{Weyl module with index -}
\def\Ind(#1,#2,#3){\mbox{Ind}_{#1}^{#2}#3} 
\gge{Ind(-1,-2,-)}{\Ind(-1,-2,-)}{induction}
\def\Res(#1,#2,#3){\mbox{Res}_{#1}^{#2}#3} 
\gge{Res(-1,-2,-)}{\Res(-1,-2,-)}{restriction}
\newcommandabbreviation{\weyl}{standard}{weyl}
\newcommandabbreviation{\mod}{\mbox{mod}}{mod}
\newcommandabbreviation{\head}{\mbox{head }}{head}
\newcommandabbreviation{\Weyl}{Weyl}{Weyl}
\def\SS(#1){{\cal S}_{#1}}                 
\gge{SS(-)}{\SS(-)}{Specht/Weyl module index -}
\def\LL(#1){{\cal L}_{#1}}                 
\gge{LL(-)}{\LL(-)}{simple module index -}
\noglossaryignore{\newline FUNCTORS/MAPS\newline}
\newcommandmacro{\Gg}{{\cal G}}{Gg}{G Functor}
\newcommandmacro{\Fg}{{\cal F}}{Fg}{F Functor}
\newcommandmacro{\ra}{\rightarrow}{ra}{}
\def\ses(#1,#2,#3){0\ra #1 \ra #2 \ra #3 \ra 0}   
\gge{ses(1,2,-)}{\ses(1,2,-)}{\hspace{.5in} short exact sequence}
\def\starr(#1){ \stackrel{ #1 }{\longrightarrow} }
\gge{starr(-)}{\starr(-)}{}
\newcommandmacro{\doublerightarrow}{\; -\!\!\! -\!\!\!\!\!\! \gg \;}
{doublerightarrow}{}
\noglossaryignore{\newline PARTITION ALGEBRA MAPS\newline}
\newcommandmacro{\smap}{s}{smap}{`inclusion' map}
\newcommandmacro{\tmap}{t}{tmap}{$ P_n -> S_n$}
\newcommandmacro{\pmap}{\psi}{pmap}{$ S_n -> P_n $}
\noglossaryignore{\newline MISC.\newline}
\def\Amap(#1){{\cal A}_{#1}}               
\gge{Amap(-)}{\Amap(-)}{}
\def\Rr(#1){R_{#1}}                        
\gge{Rr(-)}{\Rr(-)}{restriction of E}
\def\Cr(#1){C_{#1}}                        
\gge{Cr(-)}{\Cr(-)}{restriction of E to N}
\newcommandmacro{\Tm}{{\cal T}}{Tm}{Transfer Matrix}
\def\On(#1){{\cal I}_{#1}}
\gge{On(-)}{\On(-)}{}
\newcommandmacro{\UU}{\underline{\sqcup}}{UU}{}  
\newcommandmacro{\UUU}{\sqcup}{UUU}{}  
\newcommandmacro{\Vq}{V_Q^{\otimes n}}{Vq}{Potts config. space}
\def\bs(#1,#2){\mbox{{\Large $\ast$}}^{#1}_{#2}}  
\gge{bs(-,-)}{\bs(-,-)}{general plumbing multiplier}
\newcommand{\ignore}[1]{}
\gge{ignore\{-\}}{\ignore{-}}{ignore argument!}
\def\choo(#1,#2){ \left( \begin{array}{c} #1 \\ #2 \end{array} \right) } 
\gge{choo(-1,-)}{\choo(-1,-)}{choose}
\newcommand{\Qed}{$\Box$}
\gge{Qed}{\mbox{\Qed}}{QED}
\def\staq(#1){\stackrel{#1}{=}}            
\gge{staq(-)}{\staq(-)}{}
\def\stam(#1){\stackrel{#1}{\rightarrow}}  
\gge{stam(-)}{\stam(-)}{}
\def\mat{ \left( \begin{array} }    
\def\tam{ \end{array}  \right) }
\gge{mat/tam}{...}{matrix delimiters}
\newcommand{\beq}{\begin{equation} }
\def\eql(#1){ \begin{equation} \label{#1} 
}
\newcommand{\eq}{\end{equation} }
\def\eqal(#1){\begin{eqnarray} \label{#1} }
\def\eqa{\end{eqnarray} }
\def\lab(#1){\label{#1}
}
\def\prl(#1){ \begin{pr} \label{#1} 
}
\def\del(#1){ \begin{de} \label{#1} 
}
\gge{smeq\{-\}}{...}{small equation}
\gge{fneq\{-\}}{...}{very small equation}
\noglossaryignore{\newline HECKE/BLOB\newline}
\newcommandmacro{\Hnq}{H_n(q)}{Hnq}{ * freestanding symbol}
\newcommandmacro{\Hn}{H_n}{Hn}{      *-mod etc.}
\newcommandmacro{\A}{{\cal A}}{A}{}
\newcommandmacro{\Cwts}{C}{Cwts}{}
\newcommandmacro{\CA}{{\cal A}}{CA}{}

\newcommandmacro{\calA}{{\cal A}}{calA}{}
\newcommandmacro{\modi}{\mbox{Mod} }{modi}{was mod not modi!}
\newcommandmacro{\Wgen}{{\Bbb S}}{Wgen}{}
\def\ol(#1){\overline{#1}}
\newcommandmacro{\st}{\mbox{St}}{st}{}
\def\CMult(#1,#2){(#1:#2)}
\def\CM(#1,#2){( #1 : #2 )}
\def\FMult#1,#2{(#1:#2)}
\def\CF#1,#2{(#1:#2)}

\newcommandmacro{\Top}{\mbox{Top}}{Top}{}
\newcommandmacro{\Soc}{\mbox{Soc}}{Soc}{}
\newcommandmacro{\Head}{\mbox{Head}}{Head}{}
\newcommandmacro{\Filt}{{\cal F}}{Filt}{}
\newcommandmacro{\Mod}{\mbox{mod}}{Mod}{}
\newcommandmacro{\Resi}{\mbox{Res }}{Resi}{was without i!}
\newcommandmacro{\Indi}{\mbox{Ind }}{Indi}{was without i!}
\def\RR(#1,#2){R^{#1}_{#2}}   
\def\TT(#1,#2){T^{#1}_{#2}}   

\newcommandmacro{\Ann}{\mbox{Ann}}{Ann}{}
\newcommandmacro{\Cen}{\mbox{Cen}}{Cen}{}
\newcommandmacro{\End}{\mbox{End}}{End}{}
\newcommandabbreviation{\semisimple}{semisimple}{semisimple}
\newcommandabbreviation{\Bratteli}{Bratteli}{Bratteli}
\newcommandabbreviation{\JBC}{Jones Basic Construction}{JBC}
\newcommandabbreviation{\pa}{partition algebra}{pa}
\newcommandabbreviation{\TM}{transfer matrix}{TM}
\newcommandabbreviation{\PM}{Potts model}{PM}
\newcommandabbreviation{\QSC}{quantum spin chain}{QSC}
\newcommandabbreviation{\Hamiltonian}{Hamiltonian}{Hamiltonian}
\newcommandabbreviation{\YS}{Young symmetrizer}{YS}



%% file: xfig/PascalXY.eepic
\setlength{\unitlength}{0.00062500in}
\begingroup\makeatletter\ifx\SetFigFont\undefined%
\gdef\SetFigFont#1#2#3#4#5{%
  \reset@font\fontsize{#1}{#2pt}%
  \fontfamily{#3}\fontseries{#4}\fontshape{#5}%
  \selectfont}%
\fi\endgroup%
{\renewcommand{\dashlinestretch}{30}
\begin{picture}(5238,2475)(0,-10)
\put(3300,2325){\makebox(0,0)[lb]{\smash{{{\SetFigFont{9}{10.8}{\rmdefault}{\mddefault}{\updefault}1}}}}}
\put(2700,1725){\makebox(0,0)[lb]{\smash{{{\SetFigFont{9}{10.8}{\rmdefault}{\mddefault}{\updefault}1}}}}}
\put(3900,1725){\makebox(0,0)[lb]{\smash{{{\SetFigFont{9}{10.8}{\rmdefault}{\mddefault}{\updefault}1}}}}}
\put(2100,1125){\makebox(0,0)[lb]{\smash{{{\SetFigFont{9}{10.8}{\rmdefault}{\mddefault}{\updefault}1}}}}}
\put(3300,1125){\makebox(0,0)[lb]{\smash{{{\SetFigFont{9}{10.8}{\rmdefault}{\mddefault}{\updefault}2}}}}}
\put(4500,1125){\makebox(0,0)[lb]{\smash{{{\SetFigFont{9}{10.8}{\rmdefault}{\mddefault}{\updefault}1}}}}}
\path(600,2325)(600,825)
\path(570.000,945.000)(600.000,825.000)(630.000,945.000)
\path(1800,225)(4800,225)
\path(4680.000,195.000)(4800.000,225.000)(4680.000,255.000)
\put(750,2325){\makebox(0,0)[lb]{\smash{{{\SetFigFont{9}{10.8}{\rmdefault}{\mddefault}{\updefault}0}}}}}
\put(750,1725){\makebox(0,0)[lb]{\smash{{{\SetFigFont{9}{10.8}{\rmdefault}{\mddefault}{\updefault}1}}}}}
\put(750,1125){\makebox(0,0)[lb]{\smash{{{\SetFigFont{9}{10.8}{\rmdefault}{\mddefault}{\updefault}2}}}}}
\put(750,525){\makebox(0,0)[lb]{\smash{{{\SetFigFont{9}{10.8}{\rmdefault}{\mddefault}{\updefault}...}}}}}
\put(2100,0){\makebox(0,0)[lb]{\smash{{{\SetFigFont{9}{10.8}{\rmdefault}{\mddefault}{\updefault}-2}}}}}
\put(2700,0){\makebox(0,0)[lb]{\smash{{{\SetFigFont{9}{10.8}{\rmdefault}{\mddefault}{\updefault}-1}}}}}
\put(3300,0){\makebox(0,0)[lb]{\smash{{{\SetFigFont{9}{10.8}{\rmdefault}{\mddefault}{\updefault}0}}}}}
\put(3900,0){\makebox(0,0)[lb]{\smash{{{\SetFigFont{9}{10.8}{\rmdefault}{\mddefault}{\updefault}1}}}}}
\put(4500,0){\makebox(0,0)[lb]{\smash{{{\SetFigFont{9}{10.8}{\rmdefault}{\mddefault}{\updefault}2}}}}}
\put(5100,0){\makebox(0,0)[lb]{\smash{{{\SetFigFont{9}{10.8}{\rmdefault}{\mddefault}{\updefault}...}}}}}
\put(1500,0){\makebox(0,0)[lb]{\smash{{{\SetFigFont{9}{10.8}{\rmdefault}{\mddefault}{\updefault}...}}}}}
\put(3000,300){\makebox(0,0)[lb]{\smash{{{\SetFigFont{9}{10.8}{\rmdefault}{\mddefault}{\updefault}weight}}}}}
\put(0,1725){\makebox(0,0)[lb]{\smash{{{\SetFigFont{9}{10.8}{\rmdefault}{\mddefault}{\updefault}level}}}}}
\end{picture}
}

%% file: xfig/PascalBase.eepic
\setlength{\unitlength}{0.00055000in}
\begingroup\makeatletter\ifx\SetFigFont\undefined%
\gdef\SetFigFont#1#2#3#4#5{%
  \reset@font\fontsize{#1}{#2pt}%
  \fontfamily{#3}\fontseries{#4}\fontshape{#5}%
  \selectfont}%
\fi\endgroup%
{\renewcommand{\dashlinestretch}{30}
\begin{picture}(12011,14064)(0,-10)
\path(5262,12687)(6012,12687)(6012,12162)
	(5262,12162)(5262,12687)
\path(5262,11637)(6012,11637)(6012,12162)
	(5262,12162)(5262,11637)
\put(5337,12387){\makebox(0,0)[lb]{\smash{{{\SetFigFont{8}{9.6}{\rmdefault}{\mddefault}{\updefault}$U_1$}}}}}
\put(5337,11787){\makebox(0,0)[lb]{\smash{{{\SetFigFont{8}{9.6}{\rmdefault}{\mddefault}{\updefault}$eU_1$}}}}}
\path(5262,14037)(5562,14037)(5562,13662)
	(5262,13662)(5262,14037)
\put(5337,13812){\makebox(0,0)[lb]{\smash{{{\SetFigFont{8}{9.6}{\rmdefault}{\mddefault}{\updefault}1}}}}}
\path(7287,11037)(8712,11037)(8712,10062)
	(7287,10062)(7287,11037)
\path(7287,9537)(8712,9537)(8712,10062)
	(7287,10062)(7287,9537)
\put(7362,10737){\makebox(0,0)[lb]{\smash{{{\SetFigFont{8}{9.6}{\rmdefault}{\mddefault}{\updefault}$U_1eU_2U_1$}}}}}
\put(7362,9762){\makebox(0,0)[lb]{\smash{{{\SetFigFont{8}{9.6}{\rmdefault}{\mddefault}{\updefault}$eU_2U_1$}}}}}
\put(7362,10287){\makebox(0,0)[lb]{\smash{{{\SetFigFont{8}{9.6}{\rmdefault}{\mddefault}{\updefault}$eU_1eU_2U_1$}}}}}
\path(9612,12762)(9912,12762)(9912,12387)
	(9612,12387)(9612,12762)
\put(9687,12537){\makebox(0,0)[lb]{\smash{{{\SetFigFont{8}{9.6}{\rmdefault}{\mddefault}{\updefault}e}}}}}
\path(7287,13362)(7587,13362)(7587,12987)
	(7287,12987)(7287,13362)
\put(7362,13137){\makebox(0,0)[lb]{\smash{{{\SetFigFont{8}{9.6}{\rmdefault}{\mddefault}{\updefault}e}}}}}
\path(1812,12762)(2112,12762)(2112,12387)
	(1812,12387)(1812,12762)
\put(1887,12537){\makebox(0,0)[lb]{\smash{{{\SetFigFont{8}{9.6}{\rmdefault}{\mddefault}{\updefault}1}}}}}
\path(12,11037)(312,11037)(312,10662)
	(12,10662)(12,11037)
\put(87,10812){\makebox(0,0)[lb]{\smash{{{\SetFigFont{8}{9.6}{\rmdefault}{\mddefault}{\updefault}1}}}}}
\path(3537,13362)(3837,13362)(3837,12987)
	(3537,12987)(3537,13362)
\put(3612,13137){\makebox(0,0)[lb]{\smash{{{\SetFigFont{8}{9.6}{\rmdefault}{\mddefault}{\updefault}1}}}}}
\path(3537,9537)(4437,9537)(4437,10062)
	(3537,10062)(3537,9537)
\path(3537,11037)(4437,11037)(4437,10062)
	(3537,10062)(3537,11037)
\put(3612,10737){\makebox(0,0)[lb]{\smash{{{\SetFigFont{8}{9.6}{\rmdefault}{\mddefault}{\updefault}$U_1$}}}}}
\put(3612,9762){\makebox(0,0)[lb]{\smash{{{\SetFigFont{8}{9.6}{\rmdefault}{\mddefault}{\updefault}$U_2U_1$}}}}}
\put(3612,10287){\makebox(0,0)[lb]{\smash{{{\SetFigFont{8}{9.6}{\rmdefault}{\mddefault}{\updefault}$eU_1$}}}}}
\path(9612,7287)(11037,7287)(11037,7812)
	(9612,7812)(9612,7287)
\path(9612,9237)(11037,9237)(11037,7812)
	(9612,7812)(9612,9237)
\path(9537,9312)(11112,9312)(11112,7212)
	(9537,7212)(9537,9312)
\put(9687,8937){\makebox(0,0)[lb]{\smash{{{\SetFigFont{8}{9.6}{\rmdefault}{\mddefault}{\updefault}$U_1eU_2U_1$}}}}}
\put(9687,7962){\makebox(0,0)[lb]{\smash{{{\SetFigFont{8}{9.6}{\rmdefault}{\mddefault}{\updefault}$eU_2U_1$}}}}}
\put(9687,8487){\makebox(0,0)[lb]{\smash{{{\SetFigFont{8}{9.6}{\rmdefault}{\mddefault}{\updefault}$eU_1eU_2U_1$}}}}}
\put(9687,7512){\makebox(0,0)[lb]{\smash{{{\SetFigFont{8}{9.6}{\rmdefault}{\mddefault}{\updefault}$eU_3U_2U_1$}}}}}
\path(11637,11037)(11937,11037)(11937,10662)
	(11637,10662)(11637,11037)
\put(11712,10812){\makebox(0,0)[lb]{\smash{{{\SetFigFont{8}{9.6}{\rmdefault}{\mddefault}{\updefault}e}}}}}
\path(5187,7737)(6912,7737)(6912,9237)
	(5187,9237)(5187,7737)
\path(5187,6237)(6912,6237)(6912,7737)
	(5187,7737)(5187,6237)
\path(5112,9312)(6987,9312)(6987,6162)
	(5112,6162)(5112,9312)
\put(5262,8937){\makebox(0,0)[lb]{\smash{{{\SetFigFont{8}{9.6}{\rmdefault}{\mddefault}{\updefault}$U_1U_3U_2U_1$}}}}}
\put(5262,7962){\makebox(0,0)[lb]{\smash{{{\SetFigFont{8}{9.6}{\rmdefault}{\mddefault}{\updefault}$U_2U_1U_3U_2U_1$}}}}}
\put(5262,8487){\makebox(0,0)[lb]{\smash{{{\SetFigFont{8}{9.6}{\rmdefault}{\mddefault}{\updefault}$eU_1U_3U_2U_1$}}}}}
\put(5262,6462){\makebox(0,0)[lb]{\smash{{{\SetFigFont{8}{9.6}{\rmdefault}{\mddefault}{\updefault}$eU_2U_1U_3U_2U_1$}}}}}
\path(7212,9462)(6912,7212)
\put(5262,7437){\makebox(0,0)[lb]{\smash{{{\SetFigFont{8}{9.6}{\rmdefault}{\mddefault}{\updefault}$U_1eU_2U_1U_3U_2U_1$}}}}}
\put(5262,6987){\makebox(0,0)[lb]{\smash{{{\SetFigFont{8}{9.6}{\rmdefault}{\mddefault}{\updefault}$eU_1eU_2U_1U_3U_2U_1$}}}}}
\path(7287,12987)(6012,11937)
\path(5262,13662)(3837,13212)
\path(5562,13662)(7287,13212)
\path(3537,12987)(2112,12612)
\path(7587,12987)(9612,12612)
\path(5187,12762)(6087,12762)(6087,11562)
	(5187,11562)(5187,12762)
\path(1812,12387)(312,10887)
\path(9912,12387)(11637,10887)
\path(6087,11562)(7287,10512)
\path(9612,12387)(8712,9912)
\path(3837,12987)(5262,12387)
\path(5187,11562)(4437,10512)
\path(2112,12387)(3537,9837)
\path(3462,11112)(4512,11112)(4512,9462)
	(3462,9462)(3462,11112)
\path(7212,11112)(8787,11112)(8787,9462)
	(7212,9462)(7212,11112)
\path(4512,9462)(5187,8637)
\path(8787,9462)(9612,8637)
\path(11637,10662)(11037,7587)
\path(312,10662)(1287,7587)
\path(1287,7287)(2712,7287)(2712,7812)
	(1287,7812)(1287,7287)
\path(1287,9237)(2712,9237)(2712,7812)
	(1287,7812)(1287,9237)
\path(1212,9312)(2787,9312)(2787,7212)
	(1212,7212)(1212,9312)
\path(3462,9462)(2712,8637)
\path(7287,2112)(9462,2112)(9462,5187)
	(7287,5187)(7287,2112)
\path(7287,2112)(9462,2112)(9462,12)
	(7287,12)(7287,2112)
\path(6987,6162)(7287,4662)
\path(9537,7212)(9462,1587)
\put(4512,12612){\makebox(0,0)[lb]{\smash{{{\SetFigFont{8}{9.6}{\rmdefault}{\mddefault}{\updefault}$.U_1$}}}}}
\put(6087,13512){\makebox(0,0)[lb]{\smash{{{\SetFigFont{8}{9.6}{\rmdefault}{\mddefault}{\updefault}$.e$}}}}}
\put(6462,12462){\makebox(0,0)[lb]{\smash{{{\SetFigFont{8}{9.6}{\rmdefault}{\mddefault}{\updefault}$.U_1$}}}}}
\put(5937,11187){\makebox(0,0)[lb]{\smash{{{\SetFigFont{8}{9.6}{\rmdefault}{\mddefault}{\updefault}$.eU_2U_1$}}}}}
\put(8862,11337){\makebox(0,0)[lb]{\smash{{{\SetFigFont{8}{9.6}{\rmdefault}{\mddefault}{\updefault}$.U_2U_1$}}}}}
\put(2562,11337){\makebox(0,0)[lb]{\smash{{{\SetFigFont{8}{9.6}{\rmdefault}{\mddefault}{\updefault}$.U_2U_1$}}}}}
\put(6987,8937){\makebox(0,0)[lb]{\smash{{{\SetFigFont{8}{9.6}{\rmdefault}{\mddefault}{\updefault}$.U_3U_2U_1$}}}}}
\put(10812,9387){\makebox(0,0)[lb]{\smash{{{\SetFigFont{8}{9.6}{\rmdefault}{\mddefault}{\updefault}$.U_3U_2U_1$}}}}}
\put(1362,8937){\makebox(0,0)[lb]{\smash{{{\SetFigFont{8}{9.6}{\rmdefault}{\mddefault}{\updefault}$U_1$}}}}}
\put(1362,7962){\makebox(0,0)[lb]{\smash{{{\SetFigFont{8}{9.6}{\rmdefault}{\mddefault}{\updefault}$U_2U_1$}}}}}
\put(1362,8487){\makebox(0,0)[lb]{\smash{{{\SetFigFont{8}{9.6}{\rmdefault}{\mddefault}{\updefault}$eU_1$}}}}}
\put(1362,7512){\makebox(0,0)[lb]{\smash{{{\SetFigFont{8}{9.6}{\rmdefault}{\mddefault}{\updefault}$U_3U_2U_1$}}}}}
\put(4287,9162){\makebox(0,0)[lb]{\smash{{{\SetFigFont{8}{9.6}{\rmdefault}{\mddefault}{\updefault}$.U_3U_2U_1$}}}}}
\put(162,9387){\makebox(0,0)[lb]{\smash{{{\SetFigFont{8}{9.6}{\rmdefault}{\mddefault}{\updefault}$.U_3U_2U_1$}}}}}
\put(7362,3312){\makebox(0,0)[lb]{\smash{{{\SetFigFont{8}{9.6}{\rmdefault}{\mddefault}{\updefault}$U_1eU_2U_1U_3U_2U_1eU_4U_2U_1U_3$}}}}}
\put(7362,2862){\makebox(0,0)[lb]{\smash{{{\SetFigFont{8}{9.6}{\rmdefault}{\mddefault}{\updefault}$eU_1eU_2U_1U_3U_2U_1eU_4U_2U_1U_3$}}}}}
\put(7362,4812){\makebox(0,0)[lb]{\smash{{{\SetFigFont{8}{9.6}{\rmdefault}{\mddefault}{\updefault}$U_1U_3U_2U_1eU_4U_2U_1U_3$}}}}}
\put(7362,3837){\makebox(0,0)[lb]{\smash{{{\SetFigFont{8}{9.6}{\rmdefault}{\mddefault}{\updefault}$U_2U_1U_3U_2U_1eU_4U_2U_1U_3$}}}}}
\put(7362,4362){\makebox(0,0)[lb]{\smash{{{\SetFigFont{8}{9.6}{\rmdefault}{\mddefault}{\updefault}$eU_1U_3U_2U_1eU_4U_2U_1U_3$}}}}}
\put(7362,2337){\makebox(0,0)[lb]{\smash{{{\SetFigFont{8}{9.6}{\rmdefault}{\mddefault}{\updefault}$eU_2U_1U_3U_2U_1eU_4U_2U_1U_3$}}}}}
\put(7437,1737){\makebox(0,0)[lb]{\smash{{{\SetFigFont{8}{9.6}{\rmdefault}{\mddefault}{\updefault}$U_1eU_2U_1U_4U_3U_2U_1$}}}}}
\put(7437,762){\makebox(0,0)[lb]{\smash{{{\SetFigFont{8}{9.6}{\rmdefault}{\mddefault}{\updefault}$eU_2U_1U_4U_3U_2U_1$}}}}}
\put(7437,1287){\makebox(0,0)[lb]{\smash{{{\SetFigFont{8}{9.6}{\rmdefault}{\mddefault}{\updefault}$eU_1eU_2U_1U_4U_3U_2U_1$}}}}}
\put(7437,312){\makebox(0,0)[lb]{\smash{{{\SetFigFont{8}{9.6}{\rmdefault}{\mddefault}{\updefault}$eU_3U_2U_1U_4U_3U_2U_1$}}}}}
\put(6612,5562){\makebox(0,0)[lb]{\smash{{{\SetFigFont{8}{9.6}{\rmdefault}{\mddefault}{\updefault}$eU_4U_2U_1U_3$}}}}}
\put(9012,6237){\makebox(0,0)[lb]{\smash{{{\SetFigFont{8}{9.6}{\rmdefault}{\mddefault}{\updefault}$U_4U_3U_2U_1$}}}}}
\end{picture}
}

%% file: xfig/bigdiamond02.eepic
\setlength{\unitlength}{0.00020833in}
\begingroup\makeatletter\ifx\SetFigFont\undefined%
\gdef\SetFigFont#1#2#3#4#5{%
  \reset@font\fontsize{#1}{#2pt}%
  \fontfamily{#3}\fontseries{#4}\fontshape{#5}%
  \selectfont}%
\fi\endgroup%
{\renewcommand{\dashlinestretch}{30}
\begin{picture}(21377,15639)(0,-10)
\path(12,13212)(2412,15612)
\path(2412,15612)(4812,13212)
\path(12,13212)(2412,10812)(4812,13212)
\path(612,12612)(3012,15012)
\path(1212,12012)(3612,14412)
\path(1812,11412)(4212,13812)
\path(612,13812)(3012,11412)
\path(1212,14412)(3612,12012)
\path(1812,15012)(4212,12612)
\path(12,7812)(2412,10212)
\path(2412,10212)(4812,7812)
\path(12,7812)(2412,5412)(4812,7812)
\path(612,7212)(3012,9612)
\path(1212,6612)(3612,9012)
\path(1812,6012)(4212,8412)
\path(612,8412)(3012,6012)
\path(1212,9012)(3612,6612)
\path(1812,9612)(4212,7212)
\path(12,2412)(2412,4812)
\path(2412,4812)(4812,2412)
\path(12,2412)(2412,12)(4812,2412)
\path(612,1812)(3012,4212)
\path(1212,1212)(3612,3612)
\path(1812,612)(4212,3012)
\path(612,3012)(3012,612)
\path(1212,3612)(3612,1212)
\path(1812,4212)(4212,1812)
\path(11412,13212)(13812,15612)
\path(13812,15612)(16212,13212)
\path(11412,13212)(13812,10812)(16212,13212)
\path(12012,12612)(14412,15012)
\path(12612,12012)(15012,14412)
\path(13212,11412)(15612,13812)
\path(12012,13812)(14412,11412)
\path(12612,14412)(15012,12012)
\path(13212,15012)(15612,12612)
\path(11412,7812)(13812,10212)
\path(13812,10212)(16212,7812)
\path(11412,7812)(13812,5412)(16212,7812)
\path(12012,7212)(14412,9612)
\path(12612,6612)(15012,9012)
\path(13212,6012)(15612,8412)
\path(12012,8412)(14412,6012)
\path(12612,9012)(15012,6612)
\path(13212,9612)(15612,7212)
\path(11412,2412)(13812,4812)
\path(13812,4812)(16212,2412)
\path(11412,2412)(13812,12)(16212,2412)
\path(12012,1812)(14412,4212)
\path(12612,1212)(15012,3612)
\path(13212,612)(15612,3012)
\path(12012,3012)(14412,612)
\path(12612,3612)(15012,1212)
\path(13212,4212)(15612,1812)
\path(5112,13212)(7512,15612)
\path(7512,15612)(9912,13212)
\path(5112,13212)(7512,10812)(9912,13212)
\path(5712,12612)(8112,15012)
\path(6312,12012)(8712,14412)
\path(6912,11412)(9312,13812)
\path(5712,13812)(8112,11412)
\path(6312,14412)(8712,12012)
\path(6912,15012)(9312,12612)
\thicklines
\path(7512,15537)(8112,14937)(8712,14337)
	(9312,13737)(9912,13137)
\thinlines
\path(5112,2412)(7512,4812)
\path(7512,4812)(9912,2412)
\path(5112,2412)(7512,12)(9912,2412)
\path(5712,1812)(8112,4212)
\path(6312,1212)(8712,3612)
\path(6912,612)(9312,3012)
\path(5712,3012)(8112,612)
\path(6312,3612)(8712,1212)
\path(6912,4212)(9312,1812)
\thicklines
\path(7512,8937)(8112,8337)(8712,8937)
	(9312,8337)(9912,7737)
\path(7512,3537)(8712,2337)(9312,2937)(9912,2337)
\thinlines
\path(5112,7812)(7512,10212)
\path(7512,10212)(9912,7812)
\path(5112,7812)(7512,5412)(9912,7812)
\path(5712,7212)(8112,9612)
\path(6312,6612)(8712,9012)
\path(6912,6012)(9312,8412)
\path(5712,8412)(8112,6012)
\path(6312,9012)(8712,6612)
\path(6912,9612)(9312,7212)
\path(16512,13212)(18912,15612)
\path(18912,15612)(21312,13212)
\path(16512,13212)(18912,10812)(21312,13212)
\path(17112,12612)(19512,15012)
\path(17712,12012)(20112,14412)
\path(18312,11412)(20712,13812)
\path(17112,13812)(19512,11412)
\path(17712,14412)(20112,12012)
\path(18312,15012)(20712,12612)
\path(16512,7812)(18912,10212)
\path(18912,10212)(21312,7812)
\path(16512,7812)(18912,5412)(21312,7812)
\path(17112,7212)(19512,9612)
\path(17712,6612)(20112,9012)
\path(18312,6012)(20712,8412)
\path(17112,8412)(19512,6012)
\path(17712,9012)(20112,6612)
\path(18312,9612)(20712,7212)
\path(16512,2412)(18912,4812)
\path(18912,4812)(21312,2412)
\path(16512,2412)(18912,12)(21312,2412)
\path(17112,1812)(19512,4212)
\path(17712,1212)(20112,3612)
\path(18312,612)(20712,3012)
\path(17112,3012)(19512,612)
\path(17712,3612)(20112,1212)
\path(18312,4212)(20712,1812)
\thicklines
\path(21312,13137)(20712,13737)(20112,13137)
	(19512,13737)(18912,13137)
\path(21312,7737)(20712,8337)(20112,8937)
	(19512,8337)(18912,7737)
\path(18912,3537)(19512,4137)(20112,3537)
	(20712,2937)(21312,2337)
\path(2412,15537)(3012,14937)(2412,14337)
	(3012,13737)(2412,13137)
\path(2412,10137)(3012,9537)(3612,8937)
	(3012,8337)(2412,7737)
\path(2412,4737)(3012,4137)(2412,3537)
	(1812,2937)(2412,2337)
\path(13812,15537)(13212,14937)(13812,14337)
	(13212,13737)(13812,13137)
\path(13812,10137)(13212,9537)(12612,8937)
	(13212,8337)(13812,7737)
\path(13812,4737)(13212,4137)(13812,3537)
	(14412,2937)(13812,2337)
\path(3912,11712)(6012,11712)
\path(5892.000,11682.000)(6012.000,11712.000)(5892.000,11742.000)
\path(3912,6312)(6012,6312)
\path(5892.000,6282.000)(6012.000,6312.000)(5892.000,6342.000)
\path(3912,912)(6012,912)
\path(5892.000,882.000)(6012.000,912.000)(5892.000,942.000)
\path(15312,11712)(17412,11712)
\path(17292.000,11682.000)(17412.000,11712.000)(17292.000,11742.000)
\path(15312,6312)(17412,6312)
\path(17292.000,6282.000)(17412.000,6312.000)(17292.000,6342.000)
\path(15312,912)(17412,912)
\path(17292.000,882.000)(17412.000,912.000)(17292.000,942.000)
\path(3912,11637)(3912,11787)
\path(3912,6237)(3912,6387)
\path(3912,837)(3912,987)
\end{picture}
}

%% file: xfig/bigdiamond04.eepic
\setlength{\unitlength}{0.00020833in}
\begingroup\makeatletter\ifx\SetFigFont\undefined%
\gdef\SetFigFont#1#2#3#4#5{%
  \reset@font\fontsize{#1}{#2pt}%
  \fontfamily{#3}\fontseries{#4}\fontshape{#5}%
  \selectfont}%
\fi\endgroup%
{\renewcommand{\dashlinestretch}{30}
\begin{picture}(21924,21639)(0,-10)
\path(312,18912)(2712,21312)
\path(2712,21312)(5112,18912)
\path(312,18912)(2712,16512)(5112,18912)
\path(912,18312)(3312,20712)
\path(1512,17712)(3912,20112)
\path(2112,17112)(4512,19512)
\path(912,19512)(3312,17112)
\path(1512,20112)(3912,17712)
\path(2112,20712)(4512,18312)
\path(312,13512)(2712,15912)
\path(2712,15912)(5112,13512)
\path(312,13512)(2712,11112)(5112,13512)
\path(912,12912)(3312,15312)
\path(1512,12312)(3912,14712)
\path(2112,11712)(4512,14112)
\path(912,14112)(3312,11712)
\path(1512,14712)(3912,12312)
\path(2112,15312)(4512,12912)
\path(11712,18912)(14112,21312)
\path(14112,21312)(16512,18912)
\path(11712,18912)(14112,16512)(16512,18912)
\path(12312,18312)(14712,20712)
\path(12912,17712)(15312,20112)
\path(13512,17112)(15912,19512)
\path(12312,19512)(14712,17112)
\path(12912,20112)(15312,17712)
\path(13512,20712)(15912,18312)
\path(11712,13512)(14112,15912)
\path(14112,15912)(16512,13512)
\path(11712,13512)(14112,11112)(16512,13512)
\path(12312,12912)(14712,15312)
\path(12912,12312)(15312,14712)
\path(13512,11712)(15912,14112)
\path(12312,14112)(14712,11712)
\path(12912,14712)(15312,12312)
\path(13512,15312)(15912,12912)
\path(16812,18912)(19212,21312)
\path(19212,21312)(21612,18912)
\path(16812,18912)(19212,16512)(21612,18912)
\path(17412,18312)(19812,20712)
\path(18012,17712)(20412,20112)
\path(18612,17112)(21012,19512)
\path(17412,19512)(19812,17112)
\path(18012,20112)(20412,17712)
\path(18612,20712)(21012,18312)
\path(16812,13512)(19212,15912)
\path(19212,15912)(21612,13512)
\path(16812,13512)(19212,11112)(21612,13512)
\path(17412,12912)(19812,15312)
\path(18012,12312)(20412,14712)
\path(18612,11712)(21012,14112)
\path(17412,14112)(19812,11712)
\path(18012,14712)(20412,12312)
\path(18612,15312)(21012,12912)
\path(312,8112)(2712,10512)
\path(2712,10512)(5112,8112)
\path(312,8112)(2712,5712)(5112,8112)
\path(912,7512)(3312,9912)
\path(1512,6912)(3912,9312)
\path(2112,6312)(4512,8712)
\path(912,8712)(3312,6312)
\path(1512,9312)(3912,6912)
\path(2112,9912)(4512,7512)
\path(11712,8112)(14112,10512)
\path(14112,10512)(16512,8112)
\path(11712,8112)(14112,5712)(16512,8112)
\path(12312,7512)(14712,9912)
\path(12912,6912)(15312,9312)
\path(13512,6312)(15912,8712)
\path(12312,8712)(14712,6312)
\path(12912,9312)(15312,6912)
\path(13512,9912)(15912,7512)
\path(5412,8112)(7812,10512)
\path(7812,10512)(10212,8112)
\path(5412,8112)(7812,5712)(10212,8112)
\path(6012,7512)(8412,9912)
\path(6612,6912)(9012,9312)
\path(7212,6312)(9612,8712)
\path(6012,8712)(8412,6312)
\path(6612,9312)(9012,6912)
\path(7212,9912)(9612,7512)
\path(16812,8112)(19212,10512)
\path(19212,10512)(21612,8112)
\path(16812,8112)(19212,5712)(21612,8112)
\path(17412,7512)(19812,9912)
\path(18012,6912)(20412,9312)
\path(18612,6312)(21012,8712)
\path(17412,8712)(19812,6312)
\path(18012,9312)(20412,6912)
\path(18612,9912)(21012,7512)
\thicklines
\path(2712,10437)(2112,9837)(2712,9237)
	(3312,8637)(3912,8037)
\path(14112,10437)(14712,9837)(14112,9237)
	(13512,8637)(12912,8037)
\path(4212,6612)(6312,6612)
\path(6192.000,6582.000)(6312.000,6612.000)(6192.000,6642.000)
\path(15612,6612)(17712,6612)
\path(17592.000,6582.000)(17712.000,6612.000)(17592.000,6642.000)
\path(4212,6537)(4212,6687)
\path(7812,8037)(8412,8637)(9012,8037)
	(9612,7437)(10212,8037)
\path(19212,8037)(19812,7437)(20412,6837)
	(21012,7437)(21612,8037)
\thinlines
\path(312,2712)(2712,5112)
\path(2712,5112)(5112,2712)
\path(312,2712)(2712,312)(5112,2712)
\path(912,2112)(3312,4512)
\path(1512,1512)(3912,3912)
\path(2112,912)(4512,3312)
\path(912,3312)(3312,912)
\path(1512,3912)(3912,1512)
\path(2112,4512)(4512,2112)
\path(11712,2712)(14112,5112)
\path(14112,5112)(16512,2712)
\path(11712,2712)(14112,312)(16512,2712)
\path(12312,2112)(14712,4512)
\path(12912,1512)(15312,3912)
\path(13512,912)(15912,3312)
\path(12312,3312)(14712,912)
\path(12912,3912)(15312,1512)
\path(13512,4512)(15912,2112)
\path(5412,2712)(7812,5112)
\path(7812,5112)(10212,2712)
\path(5412,2712)(7812,312)(10212,2712)
\path(6012,2112)(8412,4512)
\path(6612,1512)(9012,3912)
\path(7212,912)(9612,3312)
\path(6012,3312)(8412,912)
\path(6612,3912)(9012,1512)
\path(7212,4512)(9612,2112)
\path(16812,2712)(19212,5112)
\path(19212,5112)(21612,2712)
\path(16812,2712)(19212,312)(21612,2712)
\path(17412,2112)(19812,4512)
\path(18012,1512)(20412,3912)
\path(18612,912)(21012,3312)
\path(17412,3312)(19812,912)
\path(18012,3912)(20412,1512)
\path(18612,4512)(21012,2112)
\thicklines
\path(2712,21237)(3312,20637)(2712,20037)
	(3312,19437)(3912,18837)
\path(14112,21237)(13512,20637)(14112,20037)
	(13512,19437)(12912,18837)
\path(14112,15837)(13512,15237)(12912,14637)
	(13512,14037)(12912,13437)
\path(4212,17412)(6312,17412)
\path(6192.000,17382.000)(6312.000,17412.000)(6192.000,17442.000)
\path(4212,12012)(6312,12012)
\path(6192.000,11982.000)(6312.000,12012.000)(6192.000,12042.000)
\path(15612,17412)(17712,17412)
\path(17592.000,17382.000)(17712.000,17412.000)(17592.000,17442.000)
\path(15612,12012)(17712,12012)
\path(17592.000,11982.000)(17712.000,12012.000)(17592.000,12042.000)
\path(4212,17337)(4212,17487)
\path(4212,11937)(4212,12087)
\thinlines
\path(5412,13512)(7812,15912)
\path(7812,15912)(10212,13512)
\path(5412,13512)(7812,11112)(10212,13512)
\path(6012,12912)(8412,15312)
\path(7212,11712)(9612,14112)
\path(6012,14112)(8412,11712)
\path(6612,14712)(9012,12312)
\path(5412,18912)(7812,21312)
\path(7812,21312)(10212,18912)
\path(5412,18912)(7812,16512)(10212,18912)
\path(6012,18312)(8412,20712)
\path(6612,17712)(9012,20112)
\path(7212,17112)(9612,19512)
\path(6012,19512)(8412,17112)
\path(6612,20112)(9012,17712)
\path(7212,20712)(9612,18312)
\thicklines
\path(7812,20037)(8412,19437)(9012,18837)
	(9612,18237)(10212,18837)
\path(2712,15837)(3312,15237)(3912,14637)
	(3312,14037)(3912,13437)
\thinlines
\path(6612,12312)(9012,14712)
\path(7212,15312)(9612,12912)
\thicklines
\path(7812,13437)(8412,12837)(9012,13437)
	(9612,12837)(10212,13437)
\path(21612,18837)(21012,18237)(20412,17637)
	(19812,18237)(19212,17637)
\path(21612,13437)(21012,12837)(20412,13437)
	(19812,12837)(19212,12237)
\path(2712,5037)(3312,4437)(3912,3837)
	(4512,3237)(3912,2637)
\path(14112,5037)(13512,4437)(12912,3837)
	(12312,3237)(12912,2637)
\path(4212,1212)(6312,1212)
\path(6192.000,1182.000)(6312.000,1212.000)(6192.000,1242.000)
\path(15612,1212)(17712,1212)
\path(17592.000,1182.000)(17712.000,1212.000)(17592.000,1242.000)
\path(4212,1137)(4212,1287)
\path(7812,2637)(8412,2037)(9012,2637)
	(9612,3237)(10212,2637)
\path(19212,1437)(19812,2037)(20412,2637)
	(21012,3237)(21612,2637)
\thinlines
\path(11412,21612)(21912,21612)(21912,12)
	(11412,12)(11412,21612)
\path(12,21612)(10512,21612)(10512,12)
	(12,12)(12,21612)
\end{picture}
}

%% file: xfig/bigdiamond03.eepic
\setlength{\unitlength}{0.00020833in}
\begingroup\makeatletter\ifx\SetFigFont\undefined%
\gdef\SetFigFont#1#2#3#4#5{%
  \reset@font\fontsize{#1}{#2pt}%
  \fontfamily{#3}\fontseries{#4}\fontshape{#5}%
  \selectfont}%
\fi\endgroup%
{\renewcommand{\dashlinestretch}{30}
\begin{picture}(21377,28839)(0,-10)
\path(12,8412)(2412,10812)
\path(2412,10812)(4812,8412)
\path(12,8412)(2412,6012)(4812,8412)
\path(612,7812)(3012,10212)
\path(1212,7212)(3612,9612)
\path(1812,6612)(4212,9012)
\path(612,9012)(3012,6612)
\path(1212,9612)(3612,7212)
\path(1812,10212)(4212,7812)
\path(5112,8412)(7512,10812)
\path(7512,10812)(9912,8412)
\path(5112,8412)(7512,6012)(9912,8412)
\path(5712,7812)(8112,10212)
\path(6312,7212)(8712,9612)
\path(6912,6612)(9312,9012)
\path(5712,9012)(8112,6612)
\path(6312,9612)(8712,7212)
\path(6912,10212)(9312,7812)
\thicklines
\path(2412,10737)(1812,10137)(1212,9537)
	(1812,8937)(2412,8337)
\path(3912,6912)(6012,6912)
\path(5892.000,6882.000)(6012.000,6912.000)(5892.000,6942.000)
\path(9912,8337)(9312,8937)(8712,9537)
	(8112,8937)(7512,8337)
\path(9912,8337)(9312,7737)(8712,8337)
	(8112,7737)(7512,8337)
\thinlines
\path(11412,8412)(13812,10812)
\path(13812,10812)(16212,8412)
\path(11412,8412)(13812,6012)(16212,8412)
\path(12012,7812)(14412,10212)
\path(12612,7212)(15012,9612)
\path(13212,6612)(15612,9012)
\path(12012,9012)(14412,6612)
\path(12612,9612)(15012,7212)
\path(13212,10212)(15612,7812)
\path(16512,8412)(18912,10812)
\path(18912,10812)(21312,8412)
\path(16512,8412)(18912,6012)(21312,8412)
\path(17112,7812)(19512,10212)
\path(17712,7212)(20112,9612)
\path(18312,6612)(20712,9012)
\path(17112,9012)(19512,6612)
\path(17712,9612)(20112,7212)
\path(18312,10212)(20712,7812)
\thicklines
\path(13812,10737)(14412,10137)(13812,9537)
	(13212,8937)(13812,8337)
\path(15312,6912)(17412,6912)
\path(17292.000,6882.000)(17412.000,6912.000)(17292.000,6942.000)
\path(21312,8337)(20712,7737)(20112,8337)
	(19512,7737)(18912,8337)
\path(18912,9537)(20112,8337)(20712,8937)(21312,8337)
\thinlines
\path(5712,2412)(8112,4812)
\path(8112,4812)(10512,2412)
\path(5712,2412)(8112,12)(10512,2412)
\path(6312,1812)(8712,4212)
\path(6912,1212)(9312,3612)
\path(7512,612)(9912,3012)
\path(6312,3012)(8712,612)
\path(6912,3612)(9312,1212)
\path(7512,4212)(9912,1812)
\thicklines
\path(8112,4737)(7512,4137)(8112,3537)
	(7512,2937)(8112,2337)
\path(9612,912)(11712,912)
\path(11592.000,882.000)(11712.000,912.000)(11592.000,942.000)
\path(15612,2337)(15012,2937)(14412,2337)
	(13812,2937)(13212,2337)
\thinlines
\path(10812,2412)(13212,4812)
\path(13212,4812)(15612,2412)
\path(10812,2412)(13212,12)(15612,2412)
\path(11412,1812)(13812,4212)
\path(12012,1212)(14412,3612)
\path(12612,612)(15012,3012)
\path(11412,3012)(13812,612)
\path(12012,3612)(14412,1212)
\path(12612,4212)(15012,1812)
\thicklines
\path(15612,2337)(15012,1737)(14412,2337)
	(13812,1737)(13212,2337)
\thinlines
\path(5712,26412)(8112,28812)
\path(8112,28812)(10512,26412)
\path(5712,26412)(8112,24012)(10512,26412)
\path(6312,25812)(8712,28212)
\path(6912,25212)(9312,27612)
\path(7512,24612)(9912,27012)
\path(6312,27012)(8712,24612)
\path(6912,27612)(9312,25212)
\path(7512,28212)(9912,25812)
\path(10812,26412)(13212,28812)
\path(13212,28812)(15612,26412)
\path(10812,26412)(13212,24012)(15612,26412)
\path(11412,25812)(13812,28212)
\path(12012,25212)(14412,27612)
\path(12612,24612)(15012,27012)
\path(11412,27012)(13812,24612)
\path(12012,27612)(14412,25212)
\path(12612,28212)(15012,25812)
\thicklines
\path(13212,28737)(13812,28137)(14412,27537)
	(15012,26937)(15612,26337)
\path(8112,28737)(8712,28137)(8112,27537)
	(8712,26937)(8112,26337)
\path(9612,24912)(11712,24912)
\path(11592.000,24882.000)(11712.000,24912.000)(11592.000,24942.000)
\path(15612,26337)(15012,25737)(14412,26337)
	(13812,25737)(13212,26337)
\put(13737,26412){\makebox(0,0)[lb]{\smash{{{\SetFigFont{5}{6.0}{\rmdefault}{\mddefault}{\updefault}1}}}}}
\put(14937,26412){\makebox(0,0)[lb]{\smash{{{\SetFigFont{5}{6.0}{\rmdefault}{\mddefault}{\updefault}3}}}}}
\put(14337,26937){\makebox(0,0)[lb]{\smash{{{\SetFigFont{5}{6.0}{\rmdefault}{\mddefault}{\updefault}2}}}}}
\put(13737,27537){\makebox(0,0)[lb]{\smash{{{\SetFigFont{5}{6.0}{\rmdefault}{\mddefault}{\updefault}1}}}}}
\put(13137,26937){\makebox(0,0)[lb]{\smash{{{\SetFigFont{5}{6.0}{\rmdefault}{\mddefault}{\updefault}0}}}}}
\put(13137,28212){\makebox(0,0)[lb]{\smash{{{\SetFigFont{5}{6.0}{\rmdefault}{\mddefault}{\updefault}0}}}}}
\thinlines
\path(5712,20412)(8112,22812)
\path(8112,22812)(10512,20412)
\path(5712,20412)(8112,18012)(10512,20412)
\path(6312,19812)(8712,22212)
\path(6912,19212)(9312,21612)
\path(7512,18612)(9912,21012)
\path(6312,21012)(8712,18612)
\path(6912,21612)(9312,19212)
\path(7512,22212)(9912,19812)
\path(10812,20412)(13212,22812)
\path(13212,22812)(15612,20412)
\path(10812,20412)(13212,18012)(15612,20412)
\path(11412,19812)(13812,22212)
\path(12012,19212)(14412,21612)
\path(12612,18612)(15012,21012)
\path(11412,21012)(13812,18612)
\path(12012,21612)(14412,19212)
\path(12612,22212)(15012,19812)
\thicklines
\path(8112,22737)(7512,22137)(8112,21537)
	(8712,20937)(8112,20337)
\path(9612,18912)(11712,18912)
\path(11592.000,18882.000)(11712.000,18912.000)(11592.000,18942.000)
\path(13212,21537)(13812,22137)(14412,21537)
	(15012,20937)(15612,20337)
\path(15612,20337)(15012,19737)(14412,20337)
	(13812,19737)(13212,20337)
\thinlines
\path(5712,14412)(8112,16812)
\path(8112,16812)(10512,14412)
\path(5712,14412)(8112,12012)(10512,14412)
\path(6312,13812)(8712,16212)
\path(6912,13212)(9312,15612)
\path(7512,12612)(9912,15012)
\path(6312,15012)(8712,12612)
\path(6912,15612)(9312,13212)
\path(7512,16212)(9912,13812)
\thicklines
\path(8112,16737)(8712,16137)(9312,15537)
	(8712,14937)(8112,14337)
\path(9612,12912)(11712,12912)
\path(11592.000,12882.000)(11712.000,12912.000)(11592.000,12942.000)
\path(15612,14337)(15012,13737)(14412,14337)
	(13812,13737)(13212,14337)
\path(13212,15537)(13812,14937)(14412,15537)
	(15012,14937)(15612,14337)
\thinlines
\path(10812,14412)(13212,16812)
\path(13212,16812)(15612,14412)
\path(10812,14412)(13212,12012)(15612,14412)
\path(11412,13812)(13812,16212)
\path(12012,13212)(14412,15612)
\path(12612,12612)(15012,15012)
\path(11412,15012)(13812,12612)
\path(12012,15612)(14412,13212)
\path(12612,16212)(15012,13812)
\end{picture}
}

%% file: xfig/bigdiamond05.eepic
\setlength{\unitlength}{0.00016667in}
\begingroup\makeatletter\ifx\SetFigFont\undefined%
\gdef\SetFigFont#1#2#3#4#5{%
  \reset@font\fontsize{#1}{#2pt}%
  \fontfamily{#3}\fontseries{#4}\fontshape{#5}%
  \selectfont}%
\fi\endgroup%
{\renewcommand{\dashlinestretch}{30}
\begin{picture}(28952,42039)(0,-10)
\path(7587,27612)(9387,29412)(11187,27612)
	(9387,25812)(7587,27612)
\path(7887,27312)(9687,29112)
\path(8187,27012)(9987,28812)
\path(8487,26712)(10287,28512)
\path(8787,26412)(10587,28212)
\path(9087,26112)(10887,27912)
\path(9687,26112)(7887,27912)
\path(9987,26412)(8187,28212)
\path(10287,26712)(8487,28512)
\path(10587,27012)(8787,28812)
\path(10887,27312)(9087,29112)
\path(3687,27612)(5487,29412)(7287,27612)
	(5487,25812)(3687,27612)
\path(3987,27312)(5787,29112)
\path(4287,27012)(6087,28812)
\path(4587,26712)(6387,28512)
\path(4887,26412)(6687,28212)
\path(5187,26112)(6987,27912)
\path(5787,26112)(3987,27912)
\path(6087,26412)(4287,28212)
\path(6387,26712)(4587,28512)
\path(6687,27012)(4887,28812)
\path(6987,27312)(5187,29112)
\thicklines
\path(5487,29337)(5187,29037)(4887,28737)
	(5187,28437)(5487,28137)(5787,27837)(5487,27537)
\path(9387,28137)(9687,28437)(9987,28737)
	(10287,28437)(10587,28137)(10887,27837)(11187,27537)
\path(6687,26787)(6687,26637)
\path(6687,26712)(8187,26712)
\path(8067.000,26682.000)(8187.000,26712.000)(8067.000,26742.000)
\thinlines
\path(15687,27612)(17487,29412)(19287,27612)
	(17487,25812)(15687,27612)
\path(15987,27312)(17787,29112)
\path(16287,27012)(18087,28812)
\path(16587,26712)(18387,28512)
\path(16887,26412)(18687,28212)
\path(17187,26112)(18987,27912)
\path(17787,26112)(15987,27912)
\path(18087,26412)(16287,28212)
\path(18387,26712)(16587,28512)
\path(18687,27012)(16887,28812)
\path(18987,27312)(17187,29112)
\path(11787,27612)(13587,29412)(15387,27612)
	(13587,25812)(11787,27612)
\path(12087,27312)(13887,29112)
\path(12387,27012)(14187,28812)
\path(12687,26712)(14487,28512)
\path(12987,26412)(14787,28212)
\path(13287,26112)(15087,27912)
\path(13887,26112)(12087,27912)
\path(14187,26412)(12387,28212)
\path(14487,26712)(12687,28512)
\path(14787,27012)(12987,28812)
\path(15087,27312)(13287,29112)
\thicklines
\path(13587,29337)(13887,29037)(13587,28737)
	(13287,28437)(13587,28137)(13887,27837)(13587,27537)
\path(17487,28737)(17787,28437)(18087,28137)
	(18387,28437)(18687,28137)(18987,27837)(19287,27537)
\path(14787,26787)(14787,26637)
\path(14787,26712)(16287,26712)
\path(16167.000,26682.000)(16287.000,26712.000)(16167.000,26742.000)
\thinlines
\path(11487,23412)(13287,25212)(15087,23412)
	(13287,21612)(11487,23412)
\path(11787,23112)(13587,24912)
\path(12087,22812)(13887,24612)
\path(12387,22512)(14187,24312)
\path(12687,22212)(14487,24012)
\path(12987,21912)(14787,23712)
\path(13587,21912)(11787,23712)
\path(13887,22212)(12087,24012)
\path(14187,22512)(12387,24312)
\path(14487,22812)(12687,24612)
\path(14787,23112)(12987,24912)
\path(7587,23412)(9387,25212)(11187,23412)
	(9387,21612)(7587,23412)
\path(7887,23112)(9687,24912)
\path(8187,22812)(9987,24612)
\path(8487,22512)(10287,24312)
\path(8787,22212)(10587,24012)
\path(9087,21912)(10887,23712)
\path(9687,21912)(7887,23712)
\path(9987,22212)(8187,24012)
\path(10287,22512)(8487,24312)
\path(10587,22812)(8787,24612)
\path(10887,23112)(9087,24912)
\thicklines
\path(9387,25137)(9087,24837)(9387,24537)
	(9087,24237)(9387,23937)(9687,23637)(9387,23337)
\path(13287,23937)(13587,24237)(13887,23937)
	(14187,24237)(14487,23937)(14787,23637)(15087,23337)
\path(10587,22587)(10587,22437)
\path(10587,22512)(12087,22512)
\path(11967.000,22482.000)(12087.000,22512.000)(11967.000,22542.000)
\thinlines
\path(20487,23412)(22287,25212)(24087,23412)
	(22287,21612)(20487,23412)
\path(20787,23112)(22587,24912)
\path(21087,22812)(22887,24612)
\path(21387,22512)(23187,24312)
\path(21687,22212)(23487,24012)
\path(21987,21912)(23787,23712)
\path(22587,21912)(20787,23712)
\path(22887,22212)(21087,24012)
\path(23187,22512)(21387,24312)
\path(23487,22812)(21687,24612)
\path(23787,23112)(21987,24912)
\path(16587,23412)(18387,25212)(20187,23412)
	(18387,21612)(16587,23412)
\path(16887,23112)(18687,24912)
\path(17187,22812)(18987,24612)
\path(17487,22512)(19287,24312)
\path(17787,22212)(19587,24012)
\path(18087,21912)(19887,23712)
\path(18687,21912)(16887,23712)
\path(18987,22212)(17187,24012)
\path(19287,22512)(17487,24312)
\path(19587,22812)(17787,24612)
\path(19887,23112)(18087,24912)
\thicklines
\path(18387,25137)(18687,24837)(18387,24537)
	(18687,24237)(18987,23937)(18687,23637)(18387,23337)
\path(22287,24537)(22587,24237)(22887,23937)
	(23187,23637)(23487,23937)(23787,23637)(24087,23337)
\path(19587,22587)(19587,22437)
\path(19587,22512)(21087,22512)
\path(20967.000,22482.000)(21087.000,22512.000)(20967.000,22542.000)
\thinlines
\path(7587,19212)(9387,21012)(11187,19212)
	(9387,17412)(7587,19212)
\path(7887,18912)(9687,20712)
\path(8187,18612)(9987,20412)
\path(8487,18312)(10287,20112)
\path(8787,18012)(10587,19812)
\path(9087,17712)(10887,19512)
\path(9687,17712)(7887,19512)
\path(9987,18012)(8187,19812)
\path(10287,18312)(8487,20112)
\path(10587,18612)(8787,20412)
\path(10887,18912)(9087,20712)
\path(3687,19212)(5487,21012)(7287,19212)
	(5487,17412)(3687,19212)
\path(3987,18912)(5787,20712)
\path(4287,18612)(6087,20412)
\path(4587,18312)(6387,20112)
\path(4887,18012)(6687,19812)
\path(5187,17712)(6987,19512)
\path(5787,17712)(3987,19512)
\path(6087,18012)(4287,19812)
\path(6387,18312)(4587,20112)
\path(6687,18612)(4887,20412)
\path(6987,18912)(5187,20712)
\thicklines
\path(5487,20937)(5787,20637)(6087,20337)
	(6387,20037)(6087,19737)(5787,19437)(5487,19137)
\path(9387,19737)(9687,19437)(9987,19737)
	(10287,20037)(10587,19737)(10887,19437)(11187,19137)
\path(6687,18387)(6687,18237)
\path(6687,18312)(8187,18312)
\path(8067.000,18282.000)(8187.000,18312.000)(8067.000,18342.000)
\thinlines
\path(15687,19212)(17487,21012)(19287,19212)
	(17487,17412)(15687,19212)
\path(15987,18912)(17787,20712)
\path(16287,18612)(18087,20412)
\path(16587,18312)(18387,20112)
\path(16887,18012)(18687,19812)
\path(17187,17712)(18987,19512)
\path(17787,17712)(15987,19512)
\path(18087,18012)(16287,19812)
\path(18387,18312)(16587,20112)
\path(18687,18612)(16887,20412)
\path(18987,18912)(17187,20712)
\path(11787,19212)(13587,21012)(15387,19212)
	(13587,17412)(11787,19212)
\path(12087,18912)(13887,20712)
\path(12387,18612)(14187,20412)
\path(12687,18312)(14487,20112)
\path(12987,18012)(14787,19812)
\path(13287,17712)(15087,19512)
\path(13887,17712)(12087,19512)
\path(14187,18012)(12387,19812)
\path(14487,18312)(12687,20112)
\path(14787,18612)(12987,20412)
\path(15087,18912)(13287,20712)
\thicklines
\path(13587,20937)(13287,20637)(13587,20337)
	(13887,20037)(14187,19737)(13887,19437)(13587,19137)
\path(17487,19737)(17787,20037)(18087,19737)
	(18387,19437)(18687,19737)(18987,19437)(19287,19137)
\path(14787,18387)(14787,18237)
\path(14787,18312)(16287,18312)
\path(16167.000,18282.000)(16287.000,18312.000)(16167.000,18342.000)
\thinlines
\path(7587,10812)(9387,12612)(11187,10812)
	(9387,9012)(7587,10812)
\path(7887,10512)(9687,12312)
\path(8187,10212)(9987,12012)
\path(8487,9912)(10287,11712)
\path(8787,9612)(10587,11412)
\path(9087,9312)(10887,11112)
\path(9687,9312)(7887,11112)
\path(9987,9612)(8187,11412)
\path(10287,9912)(8487,11712)
\path(10587,10212)(8787,12012)
\path(10887,10512)(9087,12312)
\path(3687,10812)(5487,12612)(7287,10812)
	(5487,9012)(3687,10812)
\path(3987,10512)(5787,12312)
\path(4287,10212)(6087,12012)
\path(4587,9912)(6387,11712)
\path(4887,9612)(6687,11412)
\path(5187,9312)(6987,11112)
\path(5787,9312)(3987,11112)
\path(6087,9612)(4287,11412)
\path(6387,9912)(4587,11712)
\path(6687,10212)(4887,12012)
\path(6987,10512)(5187,12312)
\thicklines
\path(5487,12537)(5187,12237)(4887,11937)
	(5187,11637)(4887,11337)(5187,11037)(5487,10737)
\path(9387,10737)(9687,11037)(9987,11337)
	(10287,11037)(10587,11337)(10887,11037)(11187,10737)
\path(6687,9987)(6687,9837)
\path(6687,9912)(8187,9912)
\path(8067.000,9882.000)(8187.000,9912.000)(8067.000,9942.000)
\thinlines
\path(15687,10812)(17487,12612)(19287,10812)
	(17487,9012)(15687,10812)
\path(15987,10512)(17787,12312)
\path(16287,10212)(18087,12012)
\path(16587,9912)(18387,11712)
\path(16887,9612)(18687,11412)
\path(17187,9312)(18987,11112)
\path(17787,9312)(15987,11112)
\path(18087,9612)(16287,11412)
\path(18387,9912)(16587,11712)
\path(18687,10212)(16887,12012)
\path(18987,10512)(17187,12312)
\path(11787,10812)(13587,12612)(15387,10812)
	(13587,9012)(11787,10812)
\path(12087,10512)(13887,12312)
\path(12387,10212)(14187,12012)
\path(12687,9912)(14487,11712)
\path(12987,9612)(14787,11412)
\path(13287,9312)(15087,11112)
\path(13887,9312)(12087,11112)
\path(14187,9612)(12387,11412)
\path(14487,9912)(12687,11712)
\path(14787,10212)(12987,12012)
\path(15087,10512)(13287,12312)
\thicklines
\path(13587,12537)(13887,12237)(13587,11937)
	(13287,11637)(12987,11337)(13287,11037)(13587,10737)
\path(17487,11337)(17787,11037)(18087,10737)
	(18387,11037)(18687,11337)(18987,11037)(19287,10737)
\path(14787,9987)(14787,9837)
\path(14787,9912)(16287,9912)
\path(16167.000,9882.000)(16287.000,9912.000)(16167.000,9942.000)
\thinlines
\path(23787,10812)(25587,12612)(27387,10812)
	(25587,9012)(23787,10812)
\path(24087,10512)(25887,12312)
\path(24387,10212)(26187,12012)
\path(24687,9912)(26487,11712)
\path(24987,9612)(26787,11412)
\path(25287,9312)(27087,11112)
\path(25887,9312)(24087,11112)
\path(26187,9612)(24387,11412)
\path(26487,9912)(24687,11712)
\path(26787,10212)(24987,12012)
\path(27087,10512)(25287,12312)
\path(19887,10812)(21687,12612)(23487,10812)
	(21687,9012)(19887,10812)
\path(20187,10512)(21987,12312)
\path(20487,10212)(22287,12012)
\path(20787,9912)(22587,11712)
\path(21087,9612)(22887,11412)
\path(21387,9312)(23187,11112)
\path(21987,9312)(20187,11112)
\path(22287,9612)(20487,11412)
\path(22587,9912)(20787,11712)
\path(22887,10212)(21087,12012)
\path(23187,10512)(21387,12312)
\thicklines
\path(21687,12537)(21987,12237)(22287,11937)
	(21987,11637)(21687,11337)(21387,11037)(21687,10737)
\path(25587,11337)(25887,11037)(26187,11337)
	(26487,11037)(26787,10737)(27087,11037)(27387,10737)
\path(22887,9987)(22887,9837)
\path(22887,9912)(24387,9912)
\path(24267.000,9882.000)(24387.000,9912.000)(24267.000,9942.000)
\thinlines
\path(7287,6012)(9087,7812)(10887,6012)
	(9087,4212)(7287,6012)
\path(7587,5712)(9387,7512)
\path(7887,5412)(9687,7212)
\path(8187,5112)(9987,6912)
\path(8487,4812)(10287,6612)
\path(8787,4512)(10587,6312)
\path(9387,4512)(7587,6312)
\path(9687,4812)(7887,6612)
\path(9987,5112)(8187,6912)
\path(10287,5412)(8487,7212)
\path(10587,5712)(8787,7512)
\path(3387,6012)(5187,7812)(6987,6012)
	(5187,4212)(3387,6012)
\path(3687,5712)(5487,7512)
\path(3987,5412)(5787,7212)
\path(4287,5112)(6087,6912)
\path(4587,4812)(6387,6612)
\path(4887,4512)(6687,6312)
\path(5487,4512)(3687,6312)
\path(5787,4812)(3987,6612)
\path(6087,5112)(4287,6912)
\path(6387,5412)(4587,7212)
\path(6687,5712)(4887,7512)
\thicklines
\path(5187,7737)(4887,7437)(4587,7137)
	(4887,6837)(5187,6537)(4887,6237)(5187,5937)
\path(9087,5937)(9387,6237)(9687,6537)
	(9987,6237)(10287,5937)(10587,6237)(10887,5937)
\path(6387,5187)(6387,5037)
\path(6387,5112)(7887,5112)
\path(7767.000,5082.000)(7887.000,5112.000)(7767.000,5142.000)
\thinlines
\path(15387,6012)(17187,7812)(18987,6012)
	(17187,4212)(15387,6012)
\path(15687,5712)(17487,7512)
\path(15987,5412)(17787,7212)
\path(16287,5112)(18087,6912)
\path(16587,4812)(18387,6612)
\path(16887,4512)(18687,6312)
\path(17487,4512)(15687,6312)
\path(17787,4812)(15987,6612)
\path(18087,5112)(16287,6912)
\path(18387,5412)(16587,7212)
\path(18687,5712)(16887,7512)
\path(11487,6012)(13287,7812)(15087,6012)
	(13287,4212)(11487,6012)
\path(11787,5712)(13587,7512)
\path(12087,5412)(13887,7212)
\path(12387,5112)(14187,6912)
\path(12687,4812)(14487,6612)
\path(12987,4512)(14787,6312)
\path(13587,4512)(11787,6312)
\path(13887,4812)(12087,6612)
\path(14187,5112)(12387,6912)
\path(14487,5412)(12687,7212)
\path(14787,5712)(12987,7512)
\thicklines
\path(13287,7737)(12987,7437)(13287,7137)
	(12987,6837)(12687,6537)(12987,6237)(13287,5937)
\path(17187,5937)(17487,6237)(17787,5937)
	(18087,6237)(18387,6537)(18687,6237)(18987,5937)
\path(14487,5187)(14487,5037)
\path(14487,5112)(15987,5112)
\path(15867.000,5082.000)(15987.000,5112.000)(15867.000,5142.000)
\thinlines
\path(23487,6012)(25287,7812)(27087,6012)
	(25287,4212)(23487,6012)
\path(23787,5712)(25587,7512)
\path(24087,5412)(25887,7212)
\path(24387,5112)(26187,6912)
\path(24687,4812)(26487,6612)
\path(24987,4512)(26787,6312)
\path(25587,4512)(23787,6312)
\path(25887,4812)(24087,6612)
\path(26187,5112)(24387,6912)
\path(26487,5412)(24687,7212)
\path(26787,5712)(24987,7512)
\path(19587,6012)(21387,7812)(23187,6012)
	(21387,4212)(19587,6012)
\path(19887,5712)(21687,7512)
\path(20187,5412)(21987,7212)
\path(20487,5112)(22287,6912)
\path(20787,4812)(22587,6612)
\path(21087,4512)(22887,6312)
\path(21687,4512)(19887,6312)
\path(21987,4812)(20187,6612)
\path(22287,5112)(20487,6912)
\path(22587,5412)(20787,7212)
\path(22887,5712)(21087,7512)
\thicklines
\path(21387,7737)(21687,7437)(21387,7137)
	(21087,6837)(21387,6537)(21087,6237)(21387,5937)
\path(25287,6537)(25587,6237)(25887,5937)
	(26187,6237)(26487,5937)(26787,6237)(27087,5937)
\path(22587,5187)(22587,5037)
\path(22587,5112)(24087,5112)
\path(23967.000,5082.000)(24087.000,5112.000)(23967.000,5142.000)
\thinlines
\path(15387,1812)(17187,3612)(18987,1812)
	(17187,12)(15387,1812)
\path(15687,1512)(17487,3312)
\path(15987,1212)(17787,3012)
\path(16287,912)(18087,2712)
\path(16587,612)(18387,2412)
\path(16887,312)(18687,2112)
\path(17487,312)(15687,2112)
\path(17787,612)(15987,2412)
\path(18087,912)(16287,2712)
\path(18387,1212)(16587,3012)
\path(18687,1512)(16887,3312)
\path(11487,1812)(13287,3612)(15087,1812)
	(13287,12)(11487,1812)
\path(11787,1512)(13587,3312)
\path(12087,1212)(13887,3012)
\path(12387,912)(14187,2712)
\path(12687,612)(14487,2412)
\path(12987,312)(14787,2112)
\path(13587,312)(11787,2112)
\path(13887,612)(12087,2412)
\path(14187,912)(12387,2712)
\path(14487,1212)(12687,3012)
\path(14787,1512)(12987,3312)
\thicklines
\path(13287,3537)(12987,3237)(13287,2937)
	(12987,2637)(13287,2337)(12987,2037)(13287,1737)
\path(17187,1737)(17487,2037)(17787,1737)
	(18087,2037)(18387,1737)(18687,2037)(18987,1737)
\path(14487,987)(14487,837)
\path(14487,912)(15987,912)
\path(15867.000,882.000)(15987.000,912.000)(15867.000,942.000)
\thinlines
\path(25287,19212)(27087,21012)(28887,19212)
	(27087,17412)(25287,19212)
\path(25587,18912)(27387,20712)
\path(25887,18612)(27687,20412)
\path(26187,18312)(27987,20112)
\path(26487,18012)(28287,19812)
\path(26787,17712)(28587,19512)
\path(27387,17712)(25587,19512)
\path(27687,18012)(25887,19812)
\path(27987,18312)(26187,20112)
\path(28287,18612)(26487,20412)
\path(28587,18912)(26787,20712)
\path(21387,19212)(23187,21012)(24987,19212)
	(23187,17412)(21387,19212)
\path(21687,18912)(23487,20712)
\path(21987,18612)(23787,20412)
\path(22287,18312)(24087,20112)
\path(22587,18012)(24387,19812)
\path(22887,17712)(24687,19512)
\path(23487,17712)(21687,19512)
\path(23787,18012)(21987,19812)
\path(24087,18312)(22287,20112)
\path(24387,18612)(22587,20412)
\path(24687,18912)(22887,20712)
\thicklines
\path(23187,20937)(23487,20637)(23187,20337)
	(23487,20037)(23187,19737)(22887,19437)(23187,19137)
\path(27087,20337)(27387,20037)(27687,19737)
	(27987,19437)(28287,19137)(28587,19437)(28887,19137)
\path(24387,18387)(24387,18237)
\path(24387,18312)(25887,18312)
\path(25767.000,18282.000)(25887.000,18312.000)(25767.000,18342.000)
\thinlines
\path(20187,15012)(21987,16812)(23787,15012)
	(21987,13212)(20187,15012)
\path(20487,14712)(22287,16512)
\path(20787,14412)(22587,16212)
\path(21087,14112)(22887,15912)
\path(21387,13812)(23187,15612)
\path(21687,13512)(23487,15312)
\path(22287,13512)(20487,15312)
\path(22587,13812)(20787,15612)
\path(22887,14112)(21087,15912)
\path(23187,14412)(21387,16212)
\path(23487,14712)(21687,16512)
\path(16287,15012)(18087,16812)(19887,15012)
	(18087,13212)(16287,15012)
\path(16587,14712)(18387,16512)
\path(16887,14412)(18687,16212)
\path(17187,14112)(18987,15912)
\path(17487,13812)(19287,15612)
\path(17787,13512)(19587,15312)
\path(18387,13512)(16587,15312)
\path(18687,13812)(16887,15612)
\path(18987,14112)(17187,15912)
\path(19287,14412)(17487,16212)
\path(19587,14712)(17787,16512)
\thicklines
\path(18087,16737)(17787,16437)(18087,16137)
	(18387,15837)(18087,15537)(17787,15237)(18087,14937)
\path(21987,15537)(22287,15837)(22587,15537)
	(22887,15237)(23187,14937)(23487,15237)(23787,14937)
\path(19287,14187)(19287,14037)
\path(19287,14112)(20787,14112)
\path(20667.000,14082.000)(20787.000,14112.000)(20667.000,14142.000)
\thinlines
\path(12012,15012)(13812,16812)(15612,15012)
	(13812,13212)(12012,15012)
\path(12312,14712)(14112,16512)
\path(12612,14412)(14412,16212)
\path(12912,14112)(14712,15912)
\path(13212,13812)(15012,15612)
\path(13512,13512)(15312,15312)
\path(14112,13512)(12312,15312)
\path(14412,13812)(12612,15612)
\path(14712,14112)(12912,15912)
\path(15012,14412)(13212,16212)
\path(15312,14712)(13512,16512)
\path(8112,15012)(9912,16812)(11712,15012)
	(9912,13212)(8112,15012)
\path(8412,14712)(10212,16512)
\path(8712,14412)(10512,16212)
\path(9012,14112)(10812,15912)
\path(9312,13812)(11112,15612)
\path(9612,13512)(11412,15312)
\path(10212,13512)(8412,15312)
\path(10512,13812)(8712,15612)
\path(10812,14112)(9012,15912)
\path(11112,14412)(9312,16212)
\path(11412,14712)(9612,16512)
\thicklines
\path(9912,16737)(10212,16437)(10512,16137)
	(10212,15837)(10512,15537)(10212,15237)(9912,14937)
\path(13812,15537)(14112,15237)(14412,15537)
	(14712,15237)(15012,15537)(15312,15237)(15612,14937)
\path(11112,14187)(11112,14037)
\path(11112,14112)(12612,14112)
\path(12492.000,14082.000)(12612.000,14112.000)(12492.000,14142.000)
\thinlines
\path(3912,15012)(5712,16812)(7512,15012)
	(5712,13212)(3912,15012)
\path(4212,14712)(6012,16512)
\path(4512,14412)(6312,16212)
\path(4812,14112)(6612,15912)
\path(5112,13812)(6912,15612)
\path(5412,13512)(7212,15312)
\path(6012,13512)(4212,15312)
\path(6312,13812)(4512,15612)
\path(6612,14112)(4812,15912)
\path(6912,14412)(5112,16212)
\path(7212,14712)(5412,16512)
\path(3612,15012)(1812,16812)(12,15012)
	(1812,13212)(3612,15012)
\path(3312,14712)(1512,16512)
\path(3012,14412)(1212,16212)
\path(2712,14112)(912,15912)
\path(2412,13812)(612,15612)
\path(2112,13512)(312,15312)
\path(1512,13512)(3312,15312)
\path(1212,13812)(3012,15612)
\path(912,14112)(2712,15912)
\path(612,14412)(2412,16212)
\path(312,14712)(2112,16512)
\thicklines
\path(5712,14937)(6012,15237)(6312,15537)
	(6612,15837)(6912,15537)(7212,15237)(7512,14937)
\path(3012,14187)(3012,14037)
\path(3012,14112)(4512,14112)
\path(4392.000,14082.000)(4512.000,14112.000)(4392.000,14142.000)
\path(1812,16737)(1512,16437)(1212,16137)
	(912,15837)(1212,15537)(1512,15237)(1812,14937)
\thinlines
\path(11487,31812)(13287,33612)(15087,31812)
	(13287,30012)(11487,31812)
\path(11787,31512)(13587,33312)
\path(12087,31212)(13887,33012)
\path(12387,30912)(14187,32712)
\path(12687,30612)(14487,32412)
\path(12987,30312)(14787,32112)
\path(13587,30312)(11787,32112)
\path(13887,30612)(12087,32412)
\path(14187,30912)(12387,32712)
\path(14487,31212)(12687,33012)
\path(14787,31512)(12987,33312)
\path(7587,31812)(9387,33612)(11187,31812)
	(9387,30012)(7587,31812)
\path(7887,31512)(9687,33312)
\path(8187,31212)(9987,33012)
\path(8487,30912)(10287,32712)
\path(8787,30612)(10587,32412)
\path(9087,30312)(10887,32112)
\path(9687,30312)(7887,32112)
\path(9987,30612)(8187,32412)
\path(10287,30912)(8487,32712)
\path(10587,31212)(8787,33012)
\path(10887,31512)(9087,33312)
\thicklines
\path(9387,33537)(9687,33237)(9987,32937)
	(9687,32637)(9387,32337)(9687,32037)(9387,31737)
\path(13287,32937)(13587,32637)(13887,32937)
	(14187,32637)(14487,32337)(14787,32037)(15087,31737)
\path(10587,30987)(10587,30837)
\path(10587,30912)(12087,30912)
\path(11967.000,30882.000)(12087.000,30912.000)(11967.000,30942.000)
\thinlines
\path(11487,36012)(13287,37812)(15087,36012)
	(13287,34212)(11487,36012)
\path(11787,35712)(13587,37512)
\path(12087,35412)(13887,37212)
\path(12387,35112)(14187,36912)
\path(12687,34812)(14487,36612)
\path(12987,34512)(14787,36312)
\path(13587,34512)(11787,36312)
\path(13887,34812)(12087,36612)
\path(14187,35112)(12387,36912)
\path(14487,35412)(12687,37212)
\path(14787,35712)(12987,37512)
\path(7587,36012)(9387,37812)(11187,36012)
	(9387,34212)(7587,36012)
\path(7887,35712)(9687,37512)
\path(8187,35412)(9987,37212)
\path(8487,35112)(10287,36912)
\path(8787,34812)(10587,36612)
\path(9087,34512)(10887,36312)
\path(9687,34512)(7887,36312)
\path(9987,34812)(8187,36612)
\path(10287,35112)(8487,36912)
\path(10587,35412)(8787,37212)
\path(10887,35712)(9087,37512)
\thicklines
\path(9387,37737)(9087,37437)(9387,37137)
	(9687,36837)(9387,36537)(9687,36237)(9387,35937)
\path(13287,37137)(13587,37437)(13887,37137)
	(14187,36837)(14487,36537)(14787,36237)(15087,35937)
\path(10587,35187)(10587,35037)
\path(10587,35112)(12087,35112)
\path(11967.000,35082.000)(12087.000,35112.000)(11967.000,35142.000)
\thinlines
\path(11487,40212)(13287,42012)(15087,40212)
	(13287,38412)(11487,40212)
\path(11787,39912)(13587,41712)
\path(12087,39612)(13887,41412)
\path(12387,39312)(14187,41112)
\path(12687,39012)(14487,40812)
\path(12987,38712)(14787,40512)
\path(13587,38712)(11787,40512)
\path(13887,39012)(12087,40812)
\path(14187,39312)(12387,41112)
\path(14487,39612)(12687,41412)
\path(14787,39912)(12987,41712)
\path(7587,40212)(9387,42012)(11187,40212)
	(9387,38412)(7587,40212)
\path(7887,39912)(9687,41712)
\path(8187,39612)(9987,41412)
\path(8487,39312)(10287,41112)
\path(8787,39012)(10587,40812)
\path(9087,38712)(10887,40512)
\path(9687,38712)(7887,40512)
\path(9987,39012)(8187,40812)
\path(10287,39312)(8487,41112)
\path(10587,39612)(8787,41412)
\path(10887,39912)(9087,41712)
\thicklines
\path(9387,41937)(9687,41637)(9387,41337)
	(9687,41037)(9387,40737)(9687,40437)(9387,40137)
\path(13287,41937)(13587,41637)(13887,41337)
	(14187,41037)(14487,40737)(14787,40437)(15087,40137)
\path(10587,39387)(10587,39237)
\path(10587,39312)(12087,39312)
\path(11967.000,39282.000)(12087.000,39312.000)(11967.000,39342.000)
\end{picture}
}

%% file: xfig/PascalBase02.eepic
\setlength{\unitlength}{0.00055000in}
\begingroup\makeatletter\ifx\SetFigFont\undefined%
\gdef\SetFigFont#1#2#3#4#5{%
  \reset@font\fontsize{#1}{#2pt}%
  \fontfamily{#3}\fontseries{#4}\fontshape{#5}%
  \selectfont}%
\fi\endgroup%
{\renewcommand{\dashlinestretch}{30}
\begin{picture}(12011,14139)(0,-10)
\path(5262,12762)(6012,12762)(6012,12237)
	(5262,12237)(5262,12762)
\path(5262,11712)(6012,11712)(6012,12237)
	(5262,12237)(5262,11712)
\put(5337,12462){\makebox(0,0)[lb]{\smash{{{\SetFigFont{8}{9.6}{\rmdefault}{\mddefault}{\updefault}$U_1$}}}}}
\put(5337,11862){\makebox(0,0)[lb]{\smash{{{\SetFigFont{8}{9.6}{\rmdefault}{\mddefault}{\updefault}$eU_1$}}}}}
\path(5262,14112)(5562,14112)(5562,13737)
	(5262,13737)(5262,14112)
\put(5337,13887){\makebox(0,0)[lb]{\smash{{{\SetFigFont{8}{9.6}{\rmdefault}{\mddefault}{\updefault}1}}}}}
\path(7287,11112)(8712,11112)(8712,10137)
	(7287,10137)(7287,11112)
\path(7287,9612)(8712,9612)(8712,10137)
	(7287,10137)(7287,9612)
\put(7362,10812){\makebox(0,0)[lb]{\smash{{{\SetFigFont{8}{9.6}{\rmdefault}{\mddefault}{\updefault}$U_1eU_2U_1$}}}}}
\put(7362,9837){\makebox(0,0)[lb]{\smash{{{\SetFigFont{8}{9.6}{\rmdefault}{\mddefault}{\updefault}$eU_2U_1$}}}}}
\put(7362,10362){\makebox(0,0)[lb]{\smash{{{\SetFigFont{8}{9.6}{\rmdefault}{\mddefault}{\updefault}$eU_1eU_2U_1$}}}}}
\path(9612,12837)(9912,12837)(9912,12462)
	(9612,12462)(9612,12837)
\put(9687,12612){\makebox(0,0)[lb]{\smash{{{\SetFigFont{8}{9.6}{\rmdefault}{\mddefault}{\updefault}e}}}}}
\path(7287,13437)(7587,13437)(7587,13062)
	(7287,13062)(7287,13437)
\put(7362,13212){\makebox(0,0)[lb]{\smash{{{\SetFigFont{8}{9.6}{\rmdefault}{\mddefault}{\updefault}e}}}}}
\path(1812,12837)(2112,12837)(2112,12462)
	(1812,12462)(1812,12837)
\put(1887,12612){\makebox(0,0)[lb]{\smash{{{\SetFigFont{8}{9.6}{\rmdefault}{\mddefault}{\updefault}1}}}}}
\path(12,11112)(312,11112)(312,10737)
	(12,10737)(12,11112)
\put(87,10887){\makebox(0,0)[lb]{\smash{{{\SetFigFont{8}{9.6}{\rmdefault}{\mddefault}{\updefault}1}}}}}
\path(3537,13437)(3837,13437)(3837,13062)
	(3537,13062)(3537,13437)
\put(3612,13212){\makebox(0,0)[lb]{\smash{{{\SetFigFont{8}{9.6}{\rmdefault}{\mddefault}{\updefault}1}}}}}
\path(3537,9612)(4437,9612)(4437,10137)
	(3537,10137)(3537,9612)
\path(3537,11112)(4437,11112)(4437,10137)
	(3537,10137)(3537,11112)
\put(3612,10812){\makebox(0,0)[lb]{\smash{{{\SetFigFont{8}{9.6}{\rmdefault}{\mddefault}{\updefault}$U_1$}}}}}
\put(3612,9837){\makebox(0,0)[lb]{\smash{{{\SetFigFont{8}{9.6}{\rmdefault}{\mddefault}{\updefault}$U_2U_1$}}}}}
\put(3612,10362){\makebox(0,0)[lb]{\smash{{{\SetFigFont{8}{9.6}{\rmdefault}{\mddefault}{\updefault}$eU_1$}}}}}
\path(9612,7362)(11037,7362)(11037,7887)
	(9612,7887)(9612,7362)
\path(9612,9312)(11037,9312)(11037,7887)
	(9612,7887)(9612,9312)
\path(9537,9387)(11112,9387)(11112,7287)
	(9537,7287)(9537,9387)
\put(9687,9012){\makebox(0,0)[lb]{\smash{{{\SetFigFont{8}{9.6}{\rmdefault}{\mddefault}{\updefault}$U_1eU_2U_1$}}}}}
\put(9687,8037){\makebox(0,0)[lb]{\smash{{{\SetFigFont{8}{9.6}{\rmdefault}{\mddefault}{\updefault}$eU_2U_1$}}}}}
\put(9687,8562){\makebox(0,0)[lb]{\smash{{{\SetFigFont{8}{9.6}{\rmdefault}{\mddefault}{\updefault}$eU_1eU_2U_1$}}}}}
\put(9687,7587){\makebox(0,0)[lb]{\smash{{{\SetFigFont{8}{9.6}{\rmdefault}{\mddefault}{\updefault}$eU_3U_2U_1$}}}}}
\path(11637,11112)(11937,11112)(11937,10737)
	(11637,10737)(11637,11112)
\put(11712,10887){\makebox(0,0)[lb]{\smash{{{\SetFigFont{8}{9.6}{\rmdefault}{\mddefault}{\updefault}e}}}}}
\path(7287,2187)(9462,2187)(9462,5262)
	(7287,5262)(7287,2187)
\path(7287,2187)(9462,2187)(9462,87)
	(7287,87)(7287,2187)
\put(7362,3387){\makebox(0,0)[lb]{\smash{{{\SetFigFont{8}{9.6}{\rmdefault}{\mddefault}{\updefault}$U_1eU_2U_1U_3eU_4U_2U_1U_3$}}}}}
\put(7362,2937){\makebox(0,0)[lb]{\smash{{{\SetFigFont{8}{9.6}{\rmdefault}{\mddefault}{\updefault}$eU_1eU_2U_1U_3eU_4U_2U_1U_3$}}}}}
\put(7362,4887){\makebox(0,0)[lb]{\smash{{{\SetFigFont{8}{9.6}{\rmdefault}{\mddefault}{\updefault}$U_1U_3eU_4U_2U_1U_3$}}}}}
\put(7362,3912){\makebox(0,0)[lb]{\smash{{{\SetFigFont{8}{9.6}{\rmdefault}{\mddefault}{\updefault}$U_2U_1U_3eU_4U_2U_1U_3$}}}}}
\put(7362,4437){\makebox(0,0)[lb]{\smash{{{\SetFigFont{8}{9.6}{\rmdefault}{\mddefault}{\updefault}$eU_1U_3eU_4U_2U_1U_3$}}}}}
\put(7362,2412){\makebox(0,0)[lb]{\smash{{{\SetFigFont{8}{9.6}{\rmdefault}{\mddefault}{\updefault}$eU_2U_1U_3eU_4U_2U_1U_3$}}}}}
\put(7437,1812){\makebox(0,0)[lb]{\smash{{{\SetFigFont{8}{9.6}{\rmdefault}{\mddefault}{\updefault}$U_1eU_2U_1U_4U_3$}}}}}
\put(7437,837){\makebox(0,0)[lb]{\smash{{{\SetFigFont{8}{9.6}{\rmdefault}{\mddefault}{\updefault}$eU_2U_1U_4U_3$}}}}}
\put(7437,1362){\makebox(0,0)[lb]{\smash{{{\SetFigFont{8}{9.6}{\rmdefault}{\mddefault}{\updefault}$eU_1eU_2U_1U_4U_3$}}}}}
\put(7437,387){\makebox(0,0)[lb]{\smash{{{\SetFigFont{8}{9.6}{\rmdefault}{\mddefault}{\updefault}$eU_3U_2U_1U_4U_3$}}}}}
\path(2937,2187)(4812,2187)(4812,5262)
	(2937,5262)(2937,2187)
\path(2937,2187)(4812,2187)(4812,87)
	(2937,87)(2937,2187)
\put(3012,3387){\makebox(0,0)[lb]{\smash{{{\SetFigFont{8}{9.6}{\rmdefault}{\mddefault}{\updefault}$U_1eU_2U_1U_3$}}}}}
\put(3012,2937){\makebox(0,0)[lb]{\smash{{{\SetFigFont{8}{9.6}{\rmdefault}{\mddefault}{\updefault}$eU_1eU_2U_1U_3$}}}}}
\put(3012,4887){\makebox(0,0)[lb]{\smash{{{\SetFigFont{8}{9.6}{\rmdefault}{\mddefault}{\updefault}$U_1U_3$}}}}}
\put(3012,3912){\makebox(0,0)[lb]{\smash{{{\SetFigFont{8}{9.6}{\rmdefault}{\mddefault}{\updefault}$U_2U_1U_3$}}}}}
\put(3012,4437){\makebox(0,0)[lb]{\smash{{{\SetFigFont{8}{9.6}{\rmdefault}{\mddefault}{\updefault}$eU_1U_3$}}}}}
\put(3012,2412){\makebox(0,0)[lb]{\smash{{{\SetFigFont{8}{9.6}{\rmdefault}{\mddefault}{\updefault}$eU_2U_1U_3$}}}}}
\put(3087,1812){\makebox(0,0)[lb]{\smash{{{\SetFigFont{8}{9.6}{\rmdefault}{\mddefault}{\updefault}$U_1U_4U_3$}}}}}
\put(3087,837){\makebox(0,0)[lb]{\smash{{{\SetFigFont{8}{9.6}{\rmdefault}{\mddefault}{\updefault}$U_2U_1U_4U_3$}}}}}
\put(3087,1362){\makebox(0,0)[lb]{\smash{{{\SetFigFont{8}{9.6}{\rmdefault}{\mddefault}{\updefault}$eU_1U_4U_3$}}}}}
\put(3087,387){\makebox(0,0)[lb]{\smash{{{\SetFigFont{8}{9.6}{\rmdefault}{\mddefault}{\updefault}$U_3U_2U_1U_4U_3$}}}}}
\path(7287,13062)(6012,12012)
\path(5262,13737)(3837,13287)
\path(5562,13737)(7287,13287)
\path(3537,13062)(2112,12687)
\path(7587,13062)(9612,12687)
\path(5187,12837)(6087,12837)(6087,11637)
	(5187,11637)(5187,12837)
\path(1812,12462)(312,10962)
\path(9912,12462)(11637,10962)
\path(6087,11637)(7287,10587)
\path(9612,12462)(8712,9987)
\path(3837,13062)(5262,12462)
\path(5187,11637)(4437,10587)
\path(2112,12462)(3537,9912)
\path(3462,11187)(4512,11187)(4512,9537)
	(3462,9537)(3462,11187)
\path(7212,11187)(8787,11187)(8787,9537)
	(7212,9537)(7212,11187)
\path(4512,9537)(5187,8712)
\path(8787,9537)(9612,8712)
\path(11637,10737)(11037,7662)
\path(312,10737)(1287,7662)
\path(1287,7362)(2712,7362)(2712,7887)
	(1287,7887)(1287,7362)
\path(1287,9312)(2712,9312)(2712,7887)
	(1287,7887)(1287,9312)
\path(1212,9387)(2787,9387)(2787,7287)
	(1212,7287)(1212,9387)
\path(3462,9537)(2712,8712)
\path(6987,6237)(7287,4737)
\path(9537,7287)(9462,1662)
\path(7212,9537)(6912,7287)
\path(5187,7812)(6912,7812)(6912,9312)
	(5187,9312)(5187,7812)
\path(5187,6312)(6912,6312)(6912,7812)
	(5187,7812)(5187,6312)
\path(5112,9387)(6987,9387)(6987,6237)
	(5112,6237)(5112,9387)
\path(5112,6237)(4812,4737)
\path(2787,7287)(2937,1662)
\path(2862,5337)(4887,5337)(4887,12)
	(2862,12)(2862,5337)
\path(7212,5337)(9537,5337)(9537,12)
	(7212,12)(7212,5337)
\put(4512,12687){\makebox(0,0)[lb]{\smash{{{\SetFigFont{8}{9.6}{\rmdefault}{\mddefault}{\updefault}$.U_1$}}}}}
\put(6087,13587){\makebox(0,0)[lb]{\smash{{{\SetFigFont{8}{9.6}{\rmdefault}{\mddefault}{\updefault}$.e$}}}}}
\put(6462,12537){\makebox(0,0)[lb]{\smash{{{\SetFigFont{8}{9.6}{\rmdefault}{\mddefault}{\updefault}$.U_1$}}}}}
\put(5937,11262){\makebox(0,0)[lb]{\smash{{{\SetFigFont{8}{9.6}{\rmdefault}{\mddefault}{\updefault}$.eU_2U_1$}}}}}
\put(8862,11412){\makebox(0,0)[lb]{\smash{{{\SetFigFont{8}{9.6}{\rmdefault}{\mddefault}{\updefault}$.U_2U_1$}}}}}
\put(2562,11412){\makebox(0,0)[lb]{\smash{{{\SetFigFont{8}{9.6}{\rmdefault}{\mddefault}{\updefault}$.U_2U_1$}}}}}
\put(6987,9012){\makebox(0,0)[lb]{\smash{{{\SetFigFont{8}{9.6}{\rmdefault}{\mddefault}{\updefault}$.U_3$}}}}}
\put(10812,9462){\makebox(0,0)[lb]{\smash{{{\SetFigFont{8}{9.6}{\rmdefault}{\mddefault}{\updefault}$.U_3U_2U_1$}}}}}
\put(1362,9012){\makebox(0,0)[lb]{\smash{{{\SetFigFont{8}{9.6}{\rmdefault}{\mddefault}{\updefault}$U_1$}}}}}
\put(1362,8037){\makebox(0,0)[lb]{\smash{{{\SetFigFont{8}{9.6}{\rmdefault}{\mddefault}{\updefault}$U_2U_1$}}}}}
\put(1362,8562){\makebox(0,0)[lb]{\smash{{{\SetFigFont{8}{9.6}{\rmdefault}{\mddefault}{\updefault}$eU_1$}}}}}
\put(1362,7587){\makebox(0,0)[lb]{\smash{{{\SetFigFont{8}{9.6}{\rmdefault}{\mddefault}{\updefault}$U_3U_2U_1$}}}}}
\put(4287,9237){\makebox(0,0)[lb]{\smash{{{\SetFigFont{8}{9.6}{\rmdefault}{\mddefault}{\updefault}$.U_3$}}}}}
\put(162,9462){\makebox(0,0)[lb]{\smash{{{\SetFigFont{8}{9.6}{\rmdefault}{\mddefault}{\updefault}$.U_3U_2U_1$}}}}}
\put(6612,5637){\makebox(0,0)[lb]{\smash{{{\SetFigFont{8}{9.6}{\rmdefault}{\mddefault}{\updefault}$eU_4U_2U_1U_3$}}}}}
\put(9012,6312){\makebox(0,0)[lb]{\smash{{{\SetFigFont{8}{9.6}{\rmdefault}{\mddefault}{\updefault}$U_4U_3$}}}}}
\put(5262,7512){\makebox(0,0)[lb]{\smash{{{\SetFigFont{8}{9.6}{\rmdefault}{\mddefault}{\updefault}$U_1eU_2U_1U_3$}}}}}
\put(5262,7062){\makebox(0,0)[lb]{\smash{{{\SetFigFont{8}{9.6}{\rmdefault}{\mddefault}{\updefault}$eU_1eU_2U_1U_3$}}}}}
\put(5262,9012){\makebox(0,0)[lb]{\smash{{{\SetFigFont{8}{9.6}{\rmdefault}{\mddefault}{\updefault}$U_1U_3$}}}}}
\put(5262,8037){\makebox(0,0)[lb]{\smash{{{\SetFigFont{8}{9.6}{\rmdefault}{\mddefault}{\updefault}$U_2U_1U_3$}}}}}
\put(5262,8562){\makebox(0,0)[lb]{\smash{{{\SetFigFont{8}{9.6}{\rmdefault}{\mddefault}{\updefault}$eU_1U_3$}}}}}
\put(5262,6537){\makebox(0,0)[lb]{\smash{{{\SetFigFont{8}{9.6}{\rmdefault}{\mddefault}{\updefault}$eU_2U_1U_3$}}}}}
\put(2487,6312){\makebox(0,0)[lb]{\smash{{{\SetFigFont{8}{9.6}{\rmdefault}{\mddefault}{\updefault}$U_4U_3$}}}}}
\end{picture}
}